\documentclass[12pt]{article}
\usepackage{amsmath}
\usepackage{amsfonts}
\usepackage{amssymb}
\usepackage{times}
\usepackage{epsfig}
\usepackage[english]{babel}

\newtheorem{theorem}{Theorem}[section]
\newtheorem{lemma}[theorem]{Lemma}
 \newtheorem{remark}[theorem]{Remark}
  
  \newtheorem{conjecture}[theorem]{Conjecture}
  \newtheorem{definition}[theorem]{Definition}
\newtheorem{proposition}[theorem]{Proposition}
\newtheorem{corollary}[theorem]{Corollary}
\newtheorem{example}[theorem]{Example}
\newenvironment{proof}[1][Proof]{\textbf{#1.} }{\ \rule{0.5em}{0.5em}}

\newcommand{\dist}{\mathrm{dist}}

\newcommand{\tr}{{\mathrm{tr}}}
\newcommand{\Bil}{{\mathrm{Bil}(\e)}}
 \renewcommand{\Im}{{\mathrm{Im}}}

 \renewcommand{\rho}{ {\varrho}}

 \renewcommand{\epsilon}{ {\varepsilon}}

\newcommand{\abs}[1]{\left\vert#1\right\vert}

\def\T{{\mathcal T}}
\def \e {{\varepsilon}}
\def \g {{\gamma}}
\def \ga {{\gamma}}

\def \a {{\alpha}}
\def \b {{\beta}}
\def \d {{\delta}}
\def \s {{\sigma}}
\def \la {{\lambda}}

\def \o {{\omega}}
\def \z {{\zeta}}

\def \r {{\varrho}}

\def \endp { $\hfill \square$ \vskip 8 pt}

\def \R {{\mathbb {R}}}
\def \H {{\mathbb {H}}}
\def \C {{\mathbb {C}}}

\def \N {{\mathbb N}}
\def \t {{\tau}}
\def \wt {\widetilde}

\def\HH{{\mathcal H}}

\def \p {{\partial}}

\newcommand{\barint}
{\rule[.036in]{.12in}{.009in}\kern-.16in \displaystyle\int}

\begin{document}

 \title{Stability of isometric maps in the Heisenberg group
 \footnote{MSC 22E30, 53C17, 58F10; keywords: Heisenberg group,
subRiemannian geometry, biLipschitz maps.}}

\author{Nicola Arcozzi \thanks{Partially supported by the COFIN project
``Harmonic Analysis'', funded by the Italian Minister of Research.
 }
    \, and Daniele Morbidelli}
\date{\today}
\maketitle
\centerline{\sc To the memory of Juha Heinonen}
\begin{abstract}
In this paper we prove  approximation results for biLipschitz
maps in the Heisenberg group. Namely, we show that a biLipschitz
map with biLipschitz constant close to one can be pointwise
approximated, quantitatively in any fixed ball, by an isometry.
This leads to an approximation in BMO norm for the map's Pansu
derivative. We also prove that a global quasigeodesic can be
approximated by a geodesic on any fixed segment.
\end{abstract}


\renewcommand{\theequation}{\arabic{section}.\arabic{equation}}
\section{Introduction}
  In   1961 Fritz John proved the following stability estimates.
 Let $f:\R^n\to\R^n$ be a biLipschitz map
 such that   $f(0)=0$ and  the Lipschitz constant of $f$ and $f^{-1}$ is less than
 $1+\e$, where $\e>0$ is small. Then for any ball $B= B(0,R)$,
 there is  $T\in O(n)$ such that
\begin{eqnarray}
\label{gionn}
 & |f(x)-Tx|\le C_n \e R, \quad\forall \,\, x\in B \quad\text{and}
\\&
  \label{gionn2}
  \displaystyle{ \frac{1}{\mathcal{L}(B(0,R))}\int_B|f'(x)- T|dx\le C'_n \e}.
\end{eqnarray}
Here $C_n$ and $C'_n$ are dimensional constants, $f'$ is the
differential of $f$  and $\mathcal L$  denotes the Lebesgue measure.  Estimates \eqref{gionn} and \eqref{gionn2}
and their   improvements are object of considerable interest in
geometric function theory and nonlinear elasticity; see for example \cite{Ko,R,FJM,GM,ATV,Ma,CFM}.

In this article we  study  approximation
results extending \eqref{gionn}  and \eqref{gionn2} from Euclidean space to
the Heisenberg  group  $\H =\{ (z;t)\in \C\times \R\}$ equipped
with its Lie group structure and its control distance $d$.
See Section \ref{iniziale} for all
the background.

The first issue is to establish  what the correct extensions are.
It is well known that an  isometry  $T:\H\to\H$  which
 fixes the origin has the form
$(z;t)\mapsto (Az; (\det A) t)$ where $A\in O(2)$. Moreover, a
notion of differentiability for maps in the Heisenberg group  has
been introduced by Pansu \cite{P} and the Pansu differential
 can be identified with a $2\times
 2$ matrix.
Therefore it is reasonable to guess that
 the extensions of \eqref{gionn}
 and \eqref{gionn2}  have the form
\begin{eqnarray}
\label{gionn3}
 & d\big(f(z;t),(Az, (\det A) t \big) \le C(\e) R, \quad\forall \,\, (z;t)\in B(0,R)
  \quad\text{and}
\\&
  \label{gionn4}
  \displaystyle{ \frac{1}{\mathcal{L}(B(0,R))}\int_{B(0,R)}|Jf(z;t)- A|dz dt \le C'(\e) }.
\end{eqnarray}
Here $f$ is a $(1+\e)$ biLipschitz map  from $\H$ onto itself
fixing the origin and once $R>0$ is chosen,   there is  $A\in
O(2)$ such that both estimates hold. $Jf$ is the Jacobian of $f$
in the sense of Pansu, see Section \ref{iniziale}. $d$ is the
control distance and $B(0,R)=\{ (z;t)\in \H : d((0;0),
(z;t))<R\}$.   Lebesgue measure  $\mathcal L$ is the Haar measure of $\mathbb H$.

The main goal of our paper is to prove both \eqref{gionn3} and
\eqref{gionn4}, with a {\it quantitative} estimate on the constant
$C(\e)$ and $C'(\e)$. A {\it qualitative} version of the first
inequality \eqref{gionn3} can rather easily be obtained
 by Arzel\`a's   Theorem, but it does not give
any estimate of the rate of convergence to 0 of  $C(\e)$.
  Our search for quantitative estimates
for
$C(\e)$ and $C'(\e)$,  as $\e\to 0$
  involves the understanding of a number of fine
properties of the Carnot-Carath\'eodory distance in $\H$ which may
have some independent interest in subriemannian geometry.

John's proof of   \eqref{gionn}, see \cite[Lemma
IV and Theorem 3]{J}, is rather elementary, but it heavily relies
 on the Euclidean structure on $\R^n$, in particular on the isotropic nature of its geometry.
 Due to the non isotropic structure of
the Heisenberg group, the proof of
 \eqref{gionn3} cannot be obtained so easily.
In order to get estimate \eqref{gionn3}, we  examine the behaviour
under biLipschitz maps of different subsets of $\H$ and in so
doing we consider $\H$ as a metric space, making very little use
of its differential structure.
 The 
geometry of subsets of the Heisenberg group  and  more generally
of Carnot   groups
 is very rich and intricate
and it has been object of many  recent papers. See for instance
  \cite{GHL},
  \cite{FSSC},
  \cite{BHT},    \cite{BRSC},
\cite{Ba}, \cite{AF}, just to quote a few.

We shall make a substantial use of the explicit form of the geodesics for the metric $d$;
 see Section \ref{iniziale}.   Although
their equations  are known, they are  not   easy to handle, and this
introduces several new  difficulties with respect to the Euclidean
situation.
 Geodesics in the Heisenberg group have been recently used by
several authors, in order to discuss a number of different
properties of $\H$ with its control distance. See for instance
Gaveau \cite{Gav}, Kor\'anyi \cite{Kor}, Monti and Serra Cassano
\cite{MSC},
 Ambrosio and Rigot \cite{AR},
Arcozzi and Ferrari \cite{AF}.

The first result we prove concerns  the behaviour of Heisenberg
quasigeodesics. Quasigeodesics are especially studied 
in 
hyperbolic spaces  (see, e.g. \cite{GH}, \cite{Bo}).
 It is well known that  for any
$\theta\in[0,2\pi]$, $t\in\R$, the
 straight line   $\g(s) = (se^{i\theta}; t)$, $s\in \R$,
  is a global geodesic for the control metric in
$\H,$ i.e.  $d(\g(s), \g(s'))=|s-s'|$  for $s,s'\in\R$.  Moreover,
all global geodesics have this form. A $(1+\e)-$quasigeodesic is,
by definition,  any path $\g:\R\to\H$, such that
 $(1+\e)^{-1}|s-s'|\le d(\ga(s), \ga(s'))\le (1+\e)|s-s'|$, for any
 $s,s'\in\R$.

It is known that  any quasigeodesic $\gamma$ is a horizontal path (see the definition in Section \ref{iniziale}).
Denote by $\dot \g_H=(a,b)$ the invariant components of $\dot\g$
in the standard horizontal  orthonormal frame $\{X,Y\}$:
$\dot\gamma = aX(\gamma) + bY(\gamma)$ almost everywhere. Then

\medskip\noindent{\bf Theorem A} (approximation of quasigodesics).
{\it There are $\e_0>0$ and $C_0$ absolute constants such that,
given a $(1+\e)-$quasigeodesic $\g:\R\to\H$ with $\e\le\e_0$,
 then its horizontal speed $\dot \g_H$ satisfies
\begin{equation}
\label{dos2}1-C\e^{1/2}\le  \Big|\frac{1}{\mathcal{L}(I)}\int_I \dot\gamma_H
(s)ds \Big| \le 1+ C\e,
\end{equation} }

Contrary to the Euclidean case, (\ref{dos2}) is not trivially
equivalent to the definition of quasigeodesic. A peculiarly
subriemannian consequence of (\ref{dos2}) is that, for small $\e$,  any $(1+\e)-$quasigeodesic $\gamma$ passing
through a point $P_0$ in $\mathbb H$ at time $s_0$ is forever forced to
avoid a certain metric cone; extrinsically speaking, a paraboloid, having vertex at $P_0$.
See Corollary \ref{ammira}. If   $\e=0$ in \eqref{dos2}, $\gamma$
is a global geodesic.

It is likely that the 
constant $\e^{1/2}$ does not  exhibit the right order of growth with respect to $\e$.
In proving Theorem A and all the  results stated below, we
use  the known comparison between the control distance $d$ and the
Euclidean one stated  in \eqref{altedo}, which usually is not  sharp.
 This  forces us to take several  times square roots of
 $\e$. The problem of getting sharp asymptotics as $\e\to 0$ seems
 to be rather complicated and it probably  requires
 new  ideas.

   In $\mathbb H$ there are two
different kind of Euclidean planes: laterals of two-dimensional
subgroups of $\H$ and planes with a characteristic point.
 Our
second step is studying how a biLipschitz map transforms a  plane
with a characteristic point. Up to a translation, it suffices to
consider the plane $t=0$.
 We prove the
following.

\medskip
\noindent \bf Theorem B \rm (biLipschitz image of a horizontal
plane). \it There is $\e_0>0$ and $C>0$ such that, if  $f$ is an
$(1+\e)-$biLipschitz self map of the Heisenberg group with $\e\le \e_0$ and $f(0)=0$,
 for any $R>0$
there is  $A\in O(2)$ such that
\begin{equation}
\label{expo} d(f(z;0), (Az; 0 )) \le C\e^{1/16} R  ,\quad\text {for
any $z\in\R^2$, $|z|\le R$} .\medskip
\end{equation}

\medskip

\rm Then, we examine how  a biLipschitz map transforms the $t$
axis $\{  (0; t)\in \C\times\R\}$, the center of $\H$. Recall
that, from the point of view of the metric $d$, the $t-$axis is
unrectifiable and its Hausdorff dimension is 2,
 see \eqref{assicuro}. The behaviour of the $t$-axis under
 quasiconformal mappings
 has been object of some interest. See especially Heinonen and Semmes
 \cite{HS}, Question 25. Here we show that the image of the
 $t-$axis under a $(1+\e)-$biLipschitz map  lays in a metric cone around the
 $t-$axis itself. Note that biLipschitz is a smaller class than
 quasiconformal.

\medskip \noindent\bf Theorem C \rm (biLipschitz image of the
$t-$axis). \it There exists a constant $\e_0>0$
with the following property. 
Let $f$ be a $(1+\e)-$biLipschitz map such that
$f(0)=0$ and $\e<\e_0$. Then, after possibly applying the isometry
$(x,y,t)\mapsto (x,-y,-t)$, we have, for some absolute $C>0$,
\[
d\big( f(0;t), (0;t)  \big)\le C \e^{1/32} d\big( (0;0), (0;t) \big),
\quad \forall \,\,t\in\R.
\]
\rm

Finally, combining all the results obtained, we obtain the extension   of John's pointwise
approximation theorem.

\medskip
\noindent \bf Theorem D \rm (pointwise approximation).\it
\label{Tpoint} There exist $\epsilon_0>0$ and $C>0$ such that, if
 $0<\e<\e_0$, $f$ is a $(1+\e)$-biLipschitz
 map of
$\mathbb H$, $R>0$ and $P_0$ is a fixed point in ${\mathbb H}$, then there exists
an isometry $ T $ of ${\mathbb H}$ 
 such that
\begin{equation}
\label{eqjohn} d(f(P),T(P))\le C\epsilon^{1/2^{11}} R,
\end{equation}
whenever $d(P,P_0)\le R$.
\rm

\medskip

As in the Euclidean case, Theorem D and Rademacher's Theorem, which
was proved in the
Heisenberg group by Pansu \cite{P}, imply that the Jacobian of $f$ in the sense
of Pansu
(see Section \ref{iniziale}) can be approximated by means
of an isometry.

\medskip
\noindent \bf Theorem E \rm(approximation of derivatives). \it
There are constants $\epsilon_0>0$ and $C>0$ such that, if $f$ is
$(1+\epsilon)$-biLipschitz with
$0\le\epsilon<\epsilon_0$, $f(0)=0$
 and  $Jf$ is the Jacobian matrix of $f$ in the sense of Pansu,
then for $R>0$  there exists  $A\in O(2)$ such that
\begin{equation}
\label{eqBMO} \frac{1}{\mathcal{L}(B(0,R))}\int_{B(0,R)}\|Jf(Q)-A\|dQ\le
C \epsilon^{1/2^{12}}.
\end{equation}
\rm

\medskip
Equation (\ref{eqBMO}) says that $Jf$ belongs to $BMO({\mathbb
H})$. By the John-Nirenberg inequality which holds in this setting
\cite{Bu}, local uniform 
exponential integrability
 can be easily obtained, see Corollary \ref{ircocervo}.

John's result in Euclidean space is stronger in at least two
respects. First, he only assumed $f$ to  be (locally) biLipschitz
in a bounded, open subset of ${\mathbb R}^n$. In order to avoid
further complication in the proofs, we chose to work with globally
biLipschitz maps. More importantly, John deduced the validity of
(\ref{eqjohn}) and \eqref{eqBMO} with a factor $\epsilon$ on the
right-hand side, intead of our nonsharp power of $\e$. The
example at the beginning of Section \ref{olimpia} shows
 that in $\mathbb H$ the power can not be better than
$\epsilon^{1/2}$.

John-type estimates
\eqref{gionn}
and \eqref{gionn2} in $\R^n$
 have been improved   in recent literature.
    Concerning estimate
   \eqref{gionn}, we mention
   the papers \cite{ATV},  \cite{Ma}, \cite{K}, \cite{GM}. See also the
   monograph \cite{R}. For
    estimate \eqref{gionn2} see e.g. the papers \cite{Ko},
   \cite{FJM} and \cite{CFM}.
It seems that a similar stability theory for maps in a
subriemannian settings is still essentially  lacking (with the
exception of the qualitative results in \cite{D}). Here we give a
first contribution to research in this direction.

Before closing this introduction, we mention that biLipschitz maps
are quasiconformal.  A characterization of
biLipschitz maps among quasiconformal ones has been given by Balogh,
Holopainen and Tyson \cite{BHT} by means of some  modulus
estimates.  Geometric function theory
 in homogeneous groups has been developed by several authors, see Kor\'anyi and
Reimann \cite{KR1,KR2}, Pansu \cite{P},  Heinonen and Koskela
\cite{HK}, Capogna \cite{C1,C2}, Capogna and Tang \cite{CT},
Capogna and Cowling \cite{CC}, Balogh \cite{Ba}, just to quote a
few.

The article is structured as follows. In \S2 we recall  some background
 on the geometry of  $\H$ and we prove several
 lemmata which will be
used in subsequent sections. In \S3 we  prove Theorem A. In \S 4
we prove Theorem B. In  \S5 we prove Theorem C and D. In \S6 we
prove Theorem E and in \S7 we discuss some examples. In the Appendix 
we provide an elementary proof of the known classification theorem
of $\H$'s isometries, which corresponds to the case $\e=0$ in
Theorem D. The proof of this very special case guided us towards
the proof of Theorem D and it might help the reader to follow the
general structure of the article.

\medskip After this paper was submitted, D. Isangulova and S. Vodopyanov   announced 
 that they proved by means of different techniques that Theorem D holds in the higher
 dimensional Heisenberg groups $\H^n$, $n\ge 2$. Their techniques  provide a 
sharp 
estimate of the power of $\e$ but they do not work in $\H^1$. Geometric properties of quasigeodesics and biLipschitz images of horizontal planes (see our Theorems A,B,C) do not follow from their results.

\medskip \noindent \bf Acknowledgements. \rm  We thank the referee for carefully reading  the manuscript and for several comments which improved its presentation.

\section{Preliminary facts } \label{iniziale}
 \setcounter{equation}{0}
 \noindent\bf Notation. \rm
 We write
$(x,y,t)\simeq(x+iy ; t)=(z;t)\in \R^3\simeq\C\times\R$ to denote
points in the Heisenberg group $\H$. Sometimes we use a synthetic
notation $P, Q, \dots$ to denote points in $\H$. Clearly
$O=(0,0,0)$. A map $f:\H\to \H$ will be sometimes split in its
coordinate projections as follows: $ f(z;t) = (\z(z;t); \t(z;
t))= (\xi(z;t),\eta(z;t), \t(z;t)). $

Since we will take several times the square root
of $\e>0$, we fix for brevity the notation $ \e_k =\e^{1/2^k}$, so
that $\e_{k+1}=\sqrt{\e_k}$. We denote by  $C$   positive absolute
constants. The symbol $b$ will denote any real or complex
function  bounded by an absolute constant, $|b|\le C$. Both $C$
and  $b$ may  change even in the same formula.

Finally, denote by $|v|$ the Euclidean norm of a vector
$v\in\R^n$, for $n=2,3,\dots$  Write $d_0(P) = d(0, P)$. Denote spheres by $S(P,
r)=\{ Q :d(P, Q)=r\}$.
    $S^+(0,r)=S(0,r)\cap\{ t\ge 0 \}$.   Spheres and balls centered at origin are
also  denoted by $S_r:=S(0,r)$, and $B_r:=B(0,r)$.

\medskip \noindent \bf  The control
distance in the Heisenberg group. \rm Let ${\mathbb H}={\mathbb
R}^3$ be the first Heisenberg group with the product $$
(x,y,t)\cdot(x^\prime,y^\prime,t^\prime)=(x+x^\prime,y+y^\prime,t+
t^\prime+2(x^\prime y-x y^\prime)), $$ for any
$(x,y,t),(x^\prime,y^\prime,t^\prime)\in\R^3$. Denote by $L_P$ the
left translation $L_PQ:=P\cdot Q$, $P, Q\in \H$. Consider on
${\mathbb H}$ the left invariant vector fields $
X=\partial_x+2y\partial_t$ and $ Y=\partial_y-2x\partial_t.
$
The bundle  ${\mathcal H} $  spanned by $X$ and $Y$ is called the
horizontal bundle.
 A path $\gamma:[\alpha,\beta] \to{\mathbb
H}$ is said to be be a {\it horizontal curve} if $\g$ is
absolutely continuous  and there are $a, b$ measurable functions
such that
$
\dot{\gamma}(t)=a(t)X_{\gamma(t)}+b(t)Y_{\gamma(t)}$, for a.e.
$t\in[\alpha,\beta].
$
The {\it length} of $\gamma$ is
${\mathrm{length}(\g):=\int_\alpha^\beta\sqrt{a^2(t)+b^2(t)}dt} $.
Given $P,Q\in{\mathbb H}$, the \it control distance \rm
$
d(P,Q)$ is defined as the infimum (actually minimum) among the
lengths of horizontal paths connecting $P$ and $Q$.
  Later on we will discuss the
family of the
corresponding geodesics.

The ball of center $P$ and radius $R>0$ in ${\mathbb H}$ is
denoted by $B(P,R)=\{Q\in{\mathbb H}:\ d(P,Q)<R\}$. The Lebesgue
measure $dxdydt$ on ${\mathbb H}$ is, at the same, the
bi-invariant  Haar measure on ${\mathbb H}$ and, modulo a
multiplicative constant, the Hausdorff measure $H^4_d$ associated
with $d$.    Note the exponent $4$, which comes from the
homogeneous dimension of ${\mathbb H}$. $\mathcal{L}(E)$ denotes the Lebesgue
measure of $E\subset \H$.

The control distance (see \cite{NSW})  locally satisfies the
estimates
 \begin{equation}
\label{altedo} k_1 |(z;t)- (z'; t')|\le d((z;t), (z'; t'))\le  k_2 |(z;t)- (z';
t')|^{1/2},
\end{equation}
$(z;t),(z'; t')\in K $ where $K\subset\H$ is compact  and $k_1,
k_2$ depend on $K$. More precisely,
\begin{equation}
\label{rodomonte} d((z;t), (z'; t'))\approx
|z-z^\prime|+|t'-t -2\mbox{Im}\,z\overline{z^\prime}|^{1/2},
\end{equation}
with global equivalence constants.

A map $f$ from ${\mathbb H}$ to itself is $(1+\epsilon)$-{\it
biLipschitz}, $\epsilon>0$, if
\begin{equation}
\frac{1}{1+\epsilon}d(P,Q)\le d(f(P),f(Q))\le(1+\epsilon) d(P,Q),\ P,Q\in{\mathbb H}.
\end{equation}
An isometry is a $1$-biLipschitz map from $\mathbb H$ to itself.

\medskip\bf\noindent  Isometries and dilations. \rm
 The left translations $L_P:Q\mapsto P\cdot Q$ are
isometries of the Heisenberg group and they preserve the length of
a curve. Let $\theta\in{\mathbb R}$. The rotation by an angle of
$\theta$ around the $t$-axis, is the map ${\mathcal
R}_\theta:(z;t)\mapsto(e^{i\theta}z; t)$. It is   known, see
   \cite{KR1},  \cite{C1}, \cite{T} and
\cite{Ki}, that the only isometries of ${\mathbb H}$ are the
compositions of rotations, left translations and of the map
$J:(z;t)\mapsto(\overline{z};-t)$. A simple proof of this   fact,
relying directely  on properties of geodesics, is given in the
Appendix.

The dilation with parameter $\lambda>0$ of ${\mathbb H}$   is the
map $\delta_\lambda: (z;t)\mapsto(\lambda z;\lambda^2t)$. The
length of a curve is  homogeneous of degree 1 with respect to
$\d_\lambda$, i.e.
$\mathrm{length}(\delta_\lambda(\gamma))=\lambda \,
\mathrm{length}(\gamma)$, hence the same is true for the distance
function, $d(\delta_\lambda P,\delta_\lambda Q)= \lambda d(P,Q)$.

\medskip \noindent  \bf  Pansu calculus. \rm
These notions will be used in Section \ref{mirabilandia}.   Let
$f:\H\to\H$ be a Lipschitz map. The Pansu differential $Df(P)$ of
$f$ at $P\in \H$ is the map from $\H$ to $\H$ defined by
\[
Df(P)(Q)=\lim_{ \s\to 0+ }\delta_{\s^{-1}}\big\{ f(P)^{-1}\cdot f(P\cdot \delta_\s Q)
\big\},
\]
where the limit must be  uniform in $Q$ belonging to compact
sets of $\H\simeq\R^3$. Pansu  proved that the differential of
a Lipschitz map exists  almost everywhere and it is a endomorphism of
the group
$(\H,\cdot )$ into itself which commutes with dilations.   It is rather easy to check that
  that any such  morphism of $\H$ must have the form
$
(u,v,w)\mapsto(\a u+\b v, \g u+\d v, (\a \d-\b\g )w),$ for
suitable constants $ \a,\b,\g,\d\in\R.$
Therefore it can be
identified with the matrix $A=\begin{pmatrix}
  \a & \b  \\
  \g & \d
\end{pmatrix}$ and written  as
$ (u,v,w)\mapsto \left(A\binom{u}{v}; \det(A)w\right).$ Given a
point $P$ where the differential of $f$ exists and is a group
endomorphism   which commutes  with dilations,  we  denote by $Jf(P)$ its associated   $2\times 2$
matrix, so that
\begin{equation}\label{weiss}\displaystyle{ Df(P)\begin{pmatrix}
  u \\
  v \\ w
\end{pmatrix} = \begin{pmatrix}
  Jf(P)\binom{u}{v} \\
  \det(Jf)(P) w \\
\end{pmatrix}.
}\end{equation}

Note that  a smooth function need not be differentiable in Pansu
sense: the function $f(x,y,t) = (0,0,y)$ is not differentiable at
$(0,0,0)$. Moreover the mere   existence of the Pansu differential at
a point does not ensure that the latter is a morphism of $\H$, as
the function $f(x,y,t)=(x,y,2t)$ shows at the origin.

But,
 if we know that $f$ is
differentiable   in Pansu   sense at $P$ and $Df$ is a morphism,
writing
$f(P) = (\xi(P), \eta(P), \tau(P))$,
its  Pansu Jacobian matrix has the form
\begin{equation}
\label{posticino}
Jf(P)=
 \begin{pmatrix}
   X\xi(P) & Y\xi(P) \\
   X\eta(P) & Y\eta(P) \
 \end{pmatrix}.
 \end{equation}

\medskip\noindent\bf
  Geodesics and balls. \rm
We say that a curve $\gamma:I\to{\mathbb H}$ defined on an open interval $I$ of ${\mathbb R}$,  is a {\it geodesic} if it is absolutely continuous in the Euclidean sense and for any $t\in I$ there is $J\subset I$ containing $t$ such that
for all $\alpha<\beta$, $\alpha,\beta\in{  J}$, 
$
d(\gamma(\alpha),\gamma(\beta))=
\mathrm{length}\big(\gamma|_{_{[\alpha,\beta]}}\big).
$
Let $P\in{\mathbb H}$. If $P=(z;t)$ with $z\ne0$, then there is a
unique curve $\gamma$ joining $O$ and $P$, such that
$\mathrm{length}(\gamma)=d(O,P)$. If $P=(0;t)$, there are
infinitely many curves with this property.

The explicit form of geodesics in the Heisenberg group has been
calculated by several authors: see e.g. \cite{Gav}, \cite{Kor},
\cite{Str}, \cite{Bel}, \cite{Mon}.
 For each $\phi\in{\mathbb R}$
and $\alpha\in[0,2\pi]$, we have the unit-speed geodesic from the
origin
\begin{equation} \label{pietro} \gamma_{\phi,\alpha}(s)=
\begin{cases}
x(s)=\sin(\alpha)\frac{1-\cos(\phi
  s)}{\phi}+\cos(\alpha)\frac{\sin(\phi s)}{\phi},\\
y(s)=\sin(\alpha)\frac{\sin(\phi
  s)}{\phi}-\cos(\alpha)\frac{1-\cos(\phi s)}{\phi},\\
t(s)=2\frac{\phi s-\sin(\phi s)}{\phi^2}.
\end{cases}\end{equation}
In the limiting case $\phi=0$, geodesics are straight lines. The
geodesic $\gamma_{\phi,\alpha}$ is length-minimizing for $s$
varying over any interval $I$ with $|I|\le\frac{2\pi}{|\phi|}$. We
say that   $2\pi/|\phi|$ is the {\it total lifetime}
   \rm of the geodesic $\gamma$. The geodesics between arbitrary pairs of points
   can be  obtained by
left translation. The parameter $|\phi|$ has an intrinsic
geometric meaning, because $\frac{2\pi}{|\phi|}$ is the length over
which $\gamma$ is length-minimizing. Hence, $|\phi|$ is invariant
under isometries and covariant
 under dilations. The geodesics for
which $\phi\ge0$ are the ones pointing upward (this means that as
$s$ grows, $t(s)$ grows).

From the equation of the geodesics, we obtain the equation of the
geodesic sphere  centered at the origin. We denote it by $S_r$ or
$S(0,r)$. It contains all $(z;t)$ in ${\mathbb H}$ s.t.
\begin{equation}
\label{eqball}
\begin{cases}
|z|=|z|(r,\phi)=2\frac{\sin(\phi r/2)}{\phi}\\
t=t(r,\phi)=2\frac{\phi r-\sin(\phi r)}{\phi^2}
\end{cases}
\text{or}
\begin{cases}
|z|^2=|z|^2(r,\phi)=\frac{2}{\phi^2}(1-\cos(\phi r))\\
t=t(r,\phi)=2\frac{\phi r-\sin(\phi r)}{\phi^2}
\end{cases}
\end{equation}
for some $\phi\in[-2\pi/r,2\pi/r]$. See Figure 1.
\begin{figure} [ht]
 \label{figura1}
\centerline{\includegraphics{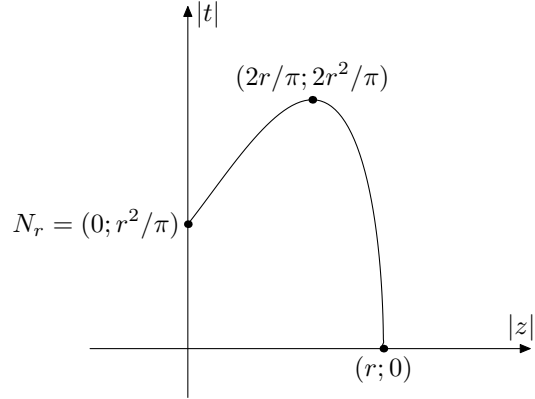} } \caption{ The sphere of
radius $r$.}
\end{figure}

\noindent Observe that
$
\phi = 0\,\Rightarrow\, (|z|;t) = (r;0),$ while $\phi =
\frac{2\pi}{r} \,\Rightarrow (|z| ; t) = (0;\frac{r^2}{\pi}),
$
so that
\begin{equation}
\label{assicuro} d((0; t), (0;s))= \sqrt{\pi \left|t-s\right|}.
\end{equation}
  The maximum and minimum values for $t$ are reached when $\phi
r=\pm\pi$ (note that $t(r,\phi)=r^2t(1,\phi r)$ and take a
derivative of $t(1,\xi)$ with respect to $\xi$).

 The upper half of the unit sphere
  $S^+(0,1):=S(0,1) \cap\{t>0\}$
 will be also written as a graph of the form
$t= u (|z|)$.
Although the function $u$ is not explicit it can be easily seen,
looking at \eqref{eqball}, that
\begin{equation} \label{assono} u(0)= \frac{1}{\pi},\quad
u\Big(\frac{2}{\pi} \Big) = \frac{2}{\pi} \quad\text{and}\quad
u'\Big(\frac{2}{\pi} \Big) = 0.
\end{equation}
A more careful look shows that $u'(0)=\frac {2}{\pi}$. The local
behaviour of $u$ near $0$ is
\begin{equation}
\label{carla} u(|z|) = \frac{1}{\pi} +\frac{2}{\pi}|z|(1 +
O(|z|)),\quad \text{  as $|z|\to 0$}.
\end{equation}
 Moreover an
easy dilation argument shows that the equation of the upper half
sphere of radius $r>0$ is
\begin{equation}
\label{ornette} \frac{t}{r^2} =u\Big(\frac{|z|}{r}\Big),\quad
|z|<r.
\end{equation}

We may also write the unit unit sphere as a graph of the $|z|$
variable, locally near $(|z|;t)=(1;0)$. The set $S(0,1)\bigcap
\Big(\C\times(-1/\pi, 1/\pi)\Big)$ can be written as $\{ (z;t):
|z| = v(t), |t|<1/ \pi \}$, where the function $v$ satisfies, for
some $C>0$,
\begin{equation}
\label{toscana}
v(|t|) = 1-C t^2 + o(t^2) , \text { as } t\to 0.
\end{equation}

A portion of the  Heisenberg ball is convex in the Euclidean sense. The following lemma is implicit in \cite{AF}, but
its proof is elementary, so we will give it here.
\begin{lemma}
\label{convesso} The convex envelope $B_{co}(O,r)$, in the
Euclidean sense, of the ball $B(O,r)$ is the solid having as
boundary the union of the portion of $S(O,r)$ corresponding to
$|\phi r|\le\pi$   in (\ref{eqball}) and the two discs $\{(z;t):\
t=\pm\frac{2}{\pi}r^2,\ |z|\le\frac{2}{\pi}r\}$.
\end{lemma}

\medskip
\noindent\it Proof. \rm   Consider the equation of
$S^+(O,r)=\partial B(O,r)\cap\{t\ge0\}$ in (\ref{eqball}). As
functions of $\phi r\in[0,2\pi]$,
  $|z|$ is  decreasing while $t$
increases on $[0,\pi]$ and decreases on $[\pi,2\pi]$. Hence, $t$
increases as $|z|$ varies in $[0,\frac{2}{\pi}r]$ and decreases as
$|z|$ varies in $[\frac{2}{\pi}r,r]$. This shows that the disc
$\{(z;t): t=\frac{2}{\pi}r^2,\ |z|\le\frac{2}{\pi}r\}$ is
contained in the convex envelope's boundary.

Let now $P=(z;t)$ be a point on $S^+(O,r)$ such that
$r\phi\in[0,\pi]$. The total lifetime of the geodesic $\gamma$
between $P$ and $O$ is $2\pi/\phi$. Since $r\phi\le \pi$, this
means that the length of the path $\gamma$ from $O$ to $P$ is less
or equal than one half of $\gamma$'s lifetime.  Consider the arc of $\gamma$ starting at $P$,
containing $O$
 and having length $\pi/\phi$, exactly one half of the lifetime of $\gamma$.
 Let $A$ be its
other endpoint. Apply now the left translation $L$ mapping $A$ to
$O$, letting $LP=P^\prime$ and $LO=O^\prime$.  Consider the ball
$LB(A, R)= B(O,R)$, where $R=\pi/\phi$. $P^\prime\in\partial
B(O,R)$ is the point in (\ref{eqball}), with $R$ instead of $r$,
corresponding to $\phi=\pi/R$: one of the point having maximum
height.  
 Hence, $B(O,R)$ stays below the set
$\{t=\frac{2}{\pi}R^2\}$, its tangent plane at $P^\prime$.
Finally, we show that $B(O^\prime,r)=L B(O,r)$ is contained in
$B(O,R)$. In fact, if $Q\in B(O^\prime,r)$, then $$ d(Q,O)\le
d(Q,O^\prime)+d(O^\prime,O)<d(P^\prime,O^\prime)+d(O^\prime,O)=d(P^\prime,O)=R,
$$ where we have used the triangle  inequality, the fact that $P^\prime\in\partial
B(O^\prime,r)$, the alignment of $P^\prime,\ O^\prime$ and $O$ on
the same length minimizing geodesic and the fact that
$P^\prime\in\partial B(O,R)$.
The inequality is strict  for any $Q\neq P'$ on the closed
ball $\{ Q: d(Q, O')\le d(P', O') \}$.
Then, $B(O^\prime,r)$ stays on one
side of its tangent plane at $P^\prime$. By translation
invariance, the same must hold with $P$ and $O$ replacing
$P^\prime$ and $O^\prime$. \endp

\medskip\noindent\bf Cones. \rm
The {\it cone} with center at $O$ and aperture $a\in{\mathbb R}$
is the set $$ \Gamma_a=\{(z;t)\in{\mathbb H}:\ t=a|z|^2\}. $$ We
could also consider the degenerate cones
$\Gamma_{\pm\infty}=\{(0;t):\ t\in{\mathbb R}^\pm\}\cup\{O\}$. The
cones centered at $O$ in ${\mathbb H}$ are the orbits of the group
generated by rotations and dilations centered at $O$, acting on
${\mathbb H}$, closed by adding the origin.
$\Gamma_{P,a}=L_P\Gamma_a$ is the cone with center at $P$ and
aperture $a$.

We introduce now a   coordinate for points in $\H$, which will be useful in Section \ref{tot}.
\begin{definition}\label{nicolino}
A point $P=(z;t)$ has coordinate $\la \ge 1$ if the geodesic
$\gamma$ starting at $O$ and passing through $P$ has a total
lifetime $\lambda\cdot d(O,P)$.
\end{definition}
Let $O^\prime$ be the other endpoint of $\gamma$, the geodesic starting at $O$ and going trough $P$. Then, $d(O,O^\prime)=\lambda\cdot d(O,P)$. The definition of $\lambda$ is
dilation invariant, $\lambda(\delta_{r}P)=\lambda(P)$, $r>0$.
 The relation between $\lambda$ and the parameter $\phi$ of $\gamma$ in
(\ref{pietro}) is
$$
\lambda=\frac{2\pi}{|\phi|r}.
$$
The points $P$ for which $\lambda(P)=\lambda$ is constant lay on the union of two cones, $\Gamma_{\pm a(\lambda)}$.
We mention that from (\ref{pietro})
or (\ref{eqball}) one deduces that $a(\lambda)\sim\frac{2\pi}{3\lambda}$ as $\lambda\to\infty$ and that
$a(\lambda)\sim\frac{1}{\pi(\lambda-1)^2}$ as $\lambda\to1$.

\begin{lemma}\label{cartozzo}  The following two facts hold.

(A)   Given $P=(z;t)\in S_1^+$, $t>0$, and $R>1$, if $\lambda(P)<R$,
then $\mathrm{dist}(P,S_R)>R-1$. Moreover, if $\la(P)\le R$, then
the distance $\mathrm{dist}(P,S_R) $   is realized by the North Pole
$N_R:=(0;R^2/\pi)\in S_R$, and by it only.    Finally, if
$\lambda(P)>R$, then $\mathrm{dist}(P,S_R)=R-1$ and the distance is realized
by a point different from the north pole.

(B) There exist $R_0>0$     and $C_0>0$, large but absolute
constants, such that, for all $R>R_0$, $(z;t)\in\H$, and
$r\in[\frac 12, 2]$,
\[
\left\{
\begin{aligned}
&\la(z;t)\ge   R
\\& (z;t)\in S_r^+
\end{aligned}  \right. \qquad \Rightarrow
\left\{
\begin{aligned}
& t \le C_0  R^{-1}
\\& 0\le r-|z| \le C_0  R^{-2}.
\end{aligned}  \right.
\]
\end{lemma}

\noindent\it Proof of (A). \rm
Suppose that $\lambda(P)<R$ and
let $Q$ be a point realizing the distance $d(Q,P)=\mathrm{dist}(P,S_R)$. 
 Then
$
1+ d(P,Q)=d(O,P)+d(P,Q)>d(O,Q)=R,
$
 the strict inequality holding
since  there is no geodesic
 passing through $O,P,Q$. Thus $d(P,Q) =\dist(P , S_R)>R-1$, as
 desired.

In order to prove the second statement,
  call ${\mathcal T}$ the closed, smooth surface obtained by taking the union of all the geodesics joining $O$ and $N_R$. Observe that
${\mathcal T}-\{N_P\}\subset B(O,R)$ and that, since $\lambda(P)\le R$, $P$
lies in the closure of the
 open set bounded by ${\mathcal T}$.
 Let $Q$ be
be any point in $S_R-\{N_R\}$. Take a geodesic $\gamma$ between
$P$ and $Q$
 and let $U$ be the the last
  point
 where $\gamma$ meets ${\mathcal T}$.
 Since $O,U$ and $N_R$ lie on the same geodesic, while
$O,U$ and $Q$ do
 not, we have that $d(O,U)+d(U,N_R)=R<d(O,U)+d(U,Q)$, hence $d(U,N_R)<d(U,Q)$. Thus,
$$
d(P,N_R)\le d(P,U)+d(U,N_R)<d(P,U)+d(U,Q)=d(P,Q).
$$
The second statement in (A) is proved.

The third statement follows easily from the definition of
lifetime of a geodesic.

\medskip \noindent \it Proof of (B). \rm   We prove (B) for $r=1$.
 The proof for $r\in[\frac 12, 2]$ is analogous.  Take a geodesic of
lifetime $R$, i.e.  with $\phi= \frac{2\pi}{R}$.
 By \eqref{eqball} with $r=1$, we have
\[
 \begin{split}
 t(1) & = \frac{R^2}{2\pi^2}\Big(  \frac{2\pi}
{R}-\sin\Big(  \frac{2\pi}{R}\Big)\Big)=\frac 23 \pi R^{-1}
+o(R^{-1}), \quad \text {while}
\\
|z(1)| & = \frac{R}{\pi}\sin\Big( \frac{\pi}{R}\Big) =
1-\frac{\pi^2}{6R^2}+o(R^{-2}).
\end{split}
\] as $R\to +\infty$.
This immediately proves the statement (B).
\endp

The properties of the set ${\mathcal T}$ in the proof are related to the fact that all points in the vertical axis are conjugate to $O$.
It was conjectured by Pansu \cite{PP,PPP}  that ${\mathcal T}$ is the extremal for the isoperimetric inequality in $\mathbb H$.

\section{Quasigeodesics in the Heisenberg group }
\label{tot} \setcounter{equation}{0}
 A (global)
$(1+\e)-$quasigeodesic in $\H$ is a curve $\ga:\R\to\H$ such that
\begin{equation}
\label{bologna}
(1+\e)^{-1}|s-\s|\le d(\ga(s), \ga(\s))\le
(1+\e)|s-\s|,\quad\text{for all} \quad s,\s\in\R.
\end{equation}
A quasigeodesic is in particular a Lipschitz embedding of $\R$ into
$\H$ equipped with the control distance. By the differentiability
theorem in \cite{P}, or by \cite[Proposition 11.4] {HK}   the path
$\g$ is horizontal, i.e. $\dot \gamma(s) = a(s)X(\gamma(s)) + b(s)
Y(\gamma(s)) $ a.e. Denote $\dot\gamma_H=(a,b)$.

\begin{theorem}
\label{bibibi}  There exist $  \e_0>0$, $C
>0$ such that for any $(1+\e)-$quasigeodesic
 $\ga$, $\e<\e_0$, and for any interval $I$ in $\mathbb R$
\begin{equation}
\label{mimma} 1-C\sqrt{\epsilon}\le\frac{1}{\mathcal{L}(I)}\left|\int_I
\dot{\gamma}_H(s)ds\right| \le 1+\epsilon.
\end{equation}
\end{theorem}
The statement of Theorem \ref{bibibi} can be explained as follows.
Let $I=[s_1,s_2]$ and $\gamma(s)=(\z(s);\t(s))$, where $\z$ is
nothing but  the Euclidean orthogonal projection  of $\gamma$ on
the plane $t=0$.
Then \eqref{mimma} reads
\begin{equation}
\label{vanga} 1-C\sqrt \e\le \frac{|\z(s_2)-\z(s_1)
|}{s_2-s_1}\le 1+\e,
\end{equation}
which implies that $\z$ is   a $(1+C\sqrt{\epsilon})-$biLipschitz
embedding of $\mathbb R$ in the plane $t=0$ endowed with the Euclidean
metric.

Theorem \ref{bibibi} clearly  will follow from the following
statement.
\begin{proposition}
\label{casimira}
 There exist $  \e_0>0$,
 $C>0$ such that for any $(1+\e)-$quasigeodesic  $\ga:\R\to \H$, $\e<\e_0$, such that $\gamma(0)=0$,
if  $\ga(s) = (\z(s); \t(s))$,
 then
\begin{equation}
\label{torino}
\left\{
\begin{aligned}
&1-C\sqrt\e\le\frac{|\z(s)|}{|s|}\le 1+C\e,
\\&|\t(s)| \le C \e^{1/4}s^2,
\end{aligned}
 \right. \quad s\in\R.
\end{equation}
\end{proposition}

We formulated here a  statement for a global geodesic.
 Actually Theorem
\ref{bibibi}   holds
 on $I=[0, L]$ for a quasigeodesic $\gamma:
\big[ -\frac{L}{\sqrt \e} , \frac{L}{\sqrt \e}
 \big]\to\H$ satisfying \eqref{bologna} for any $ s,\s$ in the
 mentioned interval.

From the previous results we get that given a
$(1+\e)-$quasigeodesic $\g=(\z;\t)$ with   $\gamma(0)=0$,
 we have for some  
$\theta\in[0,2\pi]$ and some   $b\in\R$, with $|b|\le C$,
\begin{equation}
\label{gorgo} (\z(1); \t(1))= (e^{i\theta}(1+ b \sqrt\e); b
\e^{1/4}).
\end{equation}
  Therefore, in view of \eqref{altedo}, the estimate
\begin{equation}
\label{stimantonio} d((\z(1); \t(1)), (\z(1);0))\le C\e^{1/8}
\end{equation}
holds. 
Using \eqref{stimantonio} it is easy to show that 
 the distance of $\ga(1)$ from the plane $t=0$ satisfies the  estimate
 $
\mathrm{dist}(\ga(1), \{t=0\})\le C\e^{1/4},
$
with an exponent which is better than \eqref{stimantonio}.

By dilation invariance, this last remark implies that for
$\epsilon<\epsilon_0$, a $(1+\epsilon)-$quasigeodesic $\gamma$
starting at $O$ is forced to stay outside a cone.
\begin{corollary}
\label{ammira}
There exist $  \e_0>0$, $C
>0$ such that for any $(1+\e)-$quasigeodesic  $\ga$,
 $\e<\e_0$, $\ga$ never intersects the (dilation invariant)  set $\{(z;t):
 |t| > C\e^{1/4}|z|^2\}$.
\end{corollary}

The proof of Proposition \ref{casimira} is based on the fact that
the distance between a point on the metric sphere $S(O,r)$ and the
larger, concentric sphere $S(O,R)$ can be {\it larger} than $R-r$.
  On a qualitative level, this is a
consequence of   the fact that all spheres centered at $O$ contain
points $P$ which are conjugate to $O$ along a geodesic.

\medskip
\noindent\it Proof of Proposition \ref{casimira}. \rm 
It is enough to prove the statement for $s=1$. Introduce the
numbers
\begin{equation} \label{bobo} \s= 4\sqrt\e,\quad R=
\frac{1}{\sqrt\e}.
\end{equation}
Recall that $\g(1)\in S_{\eta_1}$ and $\ga(R) \in S_{\eta_2 R }$,
where  $\eta_j\in [(1+\e)^{-1}, (1+\e)]$. Without loss of generality, one can assume that $\g(1)$ is on the northern hemisphere of $S_{\eta_1}$.   Denote by $N_{\eta_2 R}=
\Big( 0;\frac{(\eta_2 R)^2}{\pi}\Big)$ the north pole of
$S_{\eta_2 R}$. Denote $\ga(1):=(z; t)$,  recall Definition
\ref{nicolino} of the $\la-$coordinate and distinguish the
following two cases.

\smallskip\noindent \it Case A. \rm    $\la(z,t) \ge  \displaystyle{\frac{1}
{2\sqrt\e}}$.

\smallskip\noindent \it Case B. \rm    $\la(z,t) \le
\displaystyle{\frac{1} {2\sqrt\e}}$.

 In \it Case A, \rm the required estimates  \eqref{torino}
 follow immediately from
Lemma \ref{cartozzo}, part (B), which provides the estimates
$|t|\le C\sqrt \e$ and $\big|\,|z|-\eta_1\big|\le \e$ (even with
better powers than the ones in \eqref{torino}).

\smallskip The discussion of \it Case B \rm is articulated in 3
steps. The following three statements
 hold for $\e\le \e_0$, where
$\e_0$ and $C$ are absolute constants.

\smallskip\noindent \it Step B.1: \rm $ \qquad \displaystyle(z;t)
= \g(1)\in B\left( N_{\eta_2 R}, R\eta_2 - \eta_1 +\s \right) :=
B_*.
$

\smallskip \noindent\it{Step B.2:} \rm   $z$ satisfies
\begin{equation}
\label{mamma} |z| >   1-C\sqrt\e
\end{equation}
and,   as a consequence,  the first line of (\ref{torino}) holds.

\smallskip\noindent \it Step B.3: \rm  $t$ satisfies
 the estimate $ |t|\le C \e^{1/4},
$
so that the second line
of \eqref{torino} holds too.

\begin{figure} [ht]
 \centerline{\includegraphics{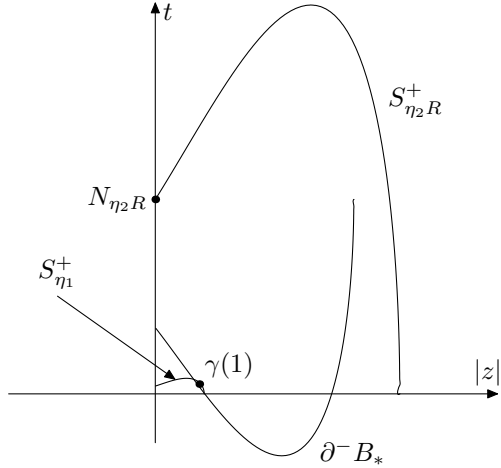} } \caption{ A graphic description of {\it Case B.}}
\end{figure}

\noindent\it Proof of Step B.1. \rm Assume by contradiction
 that $\ga(1)=(z;t)\notin
B_*$. Then it would be
$
d\left(\g(1), N_{\eta_2 R} \right) >
R\eta_2 -\eta_1 +\s.
$
Since we are assuming $ \lambda(z;t) \le \frac{1}{2\sqrt \e}
=\frac 12 R$, if $\e$ is small enough we can assert by Lemma
\ref{cartozzo}, part A, that the distance between $\g(1)$ and
$S_{\eta_2 R} $ is realized by  the point $N_{\eta_2 R} $.
 Then
\[
d(\g(R), \g(1))\ge \dist (\g(1), S_{\eta_2 R})
=d\left(\g(1), N_{\eta_2 R}
\right) > R\eta_2 -\eta_1
+\s ,
\]
because we are assuming $\g(1)\notin B_*$.
On the other hand,  the quasigeodesic
property gives $d(\g(R), \g(1))\le (R-1)(1+\e)$. Thus we get
$
 R\eta_2 -\eta_1
+\s \le (R-1)(1+\e).
$
Since
$
\eta_j\in[(1+\e)^{-1}, 1+\e]$ and $R$ and $\s$ are prescribed in
\eqref{bobo}, we get
\[
\begin{split}
\frac{1}{\sqrt \e}\eta_2    & - \eta_1 + 4\sqrt \e \le \Big(
\frac{1}{\sqrt\e} -1\Big)(1+\e)
\\&
   \Rightarrow \frac{1}{(1+\e)\sqrt\e} - (1+\e) + 4\sqrt\e\le  \Big(
\frac{1}{\sqrt\e} -1\Big)(1+\e).
\end{split}
\]
It is now  easy  to see that the last inequality can not hold for
small $\e$. This contradiction  finishes the proof of \it{Step
B.1}. \rm

\smallskip \noindent \it Proof of Step B.2. \rm
 The idea here   is to study the shape of the boundary of the
  ball $B_*=
 B \left( \Big(0; \frac{(\eta_2 R)^2}{\pi}  \Big), R\eta_2 - \eta_1 +\s
\right) $, for small  $\e$, with $R$ and $\s$ given by
\eqref{bobo}. Recall that the center of $B_*$ is $N_{\eta_2 R}$.
 Note that if we would choose $\s=0$, then the
boundary of $B_*$ would touch $\p B(0, \eta_1)$ exactly along a
circle. Choosing $\s>0$ we enlarge the ball exactly by  the amount prescribed by
 \eqref{bobo}. The next computation provides some
information about the intersection of the two balls.
In Figure 2, we represented the small upper hemisphere
$S_{\eta_1}^+$, 
the hemisphere $\p^- B_*$, the lower boundary of $B_*$ and the largest upper
 hemisphere $S_{\eta_2R}^+$.  The point $(z;t)$ belongs to the very small
region given by the intersection of $B_*$ and $B_{\eta_1}$.

The equation of the (lower hemi)sphere $S^-(0,\r)$ is  $t=  -\r^2
u\Big( \frac{z}{\r}\Big),$ $|z|\le \r$, where $u$ satisfies
\eqref{carla}.
 Taking
$\r= \eta_2 R -\eta_1+\s$ and translating upwards by the amount
$\frac{(\eta_2 R)^2}{\pi}$, we get that any point belonging to
$B_*$ should satisfy
\begin{equation}
\label{pin}
\begin{split}
t
& > \frac{(\eta_2 R)^2}{\pi} - (\eta_2 R - \eta_1 +\s)^2 u\left(
\frac{|z|}{\eta_2 R -\eta_1 +\s}\right)
\\&
 =
   \frac{(\eta_2 R)^2}{\pi} - (\eta_2 R - \eta_1 +\s)^2 \left[
  \frac{1}{\pi} +\frac{2}{\pi}\, \frac{|z|}{\eta_2 R -\eta_1 +\s}
   \Big( 1+O\Big( \frac{1}{R} \Big)\Big)\right]
\\&
    = \frac{1}{\pi} \left[2(\eta_1 -\s)\eta_2 R - (\eta_1-\s)^2 \right]
    - \frac{2}{\pi}\big(\eta_2 R - \eta_1 +\s\big)\, (1+O(1/R)) \,\,
    |z|.
\end{split}
\end{equation}
We have used here expansion
 \eqref{carla}, then we made  only algebraic simplifications and we wrote
$O(1/R)$ instead of $O(|z|/ (\eta_2 R-\eta_1 + \s))$ (this is
correct because we know that $|z|\le 2$ and we may choose  $R=
\e^{-1/2}$ large enough).

To prove \eqref{mamma}, use the fact
 that $\eta_j\in[(1+\e)^{-1}, (1+\e)].$
Thus \eqref{pin} implies
\[
\begin{split}
t
&>   \frac{1}{\pi} \left[2\big((1+\e)^{-1} -4\sqrt\e\big)
  \frac{1}{(1+\e)\sqrt\e} -
  (1+\e-4\sqrt\e)^2 \right]
 \\&
 \qquad   - \frac{2}{\pi}\Big(\frac{1+\e}{\sqrt\e}  - \frac{1}{1+\e}
  +4\sqrt\e \Big)\, (1+O(\sqrt \e)) \,\,
    |z|  .
\end{split}
\]
Some short computations show that, as $\e\to 0$,
\[
\begin{split}
2\big((1+\e)^{-1} -4\sqrt\e\big)
  \frac{1}{(1+\e)\sqrt\e} & =\frac{2}{\sqrt\e} (1-4\sqrt\e
  +O(\e)),
\end{split}
\]
and $ (1+\e-4\sqrt\e)^2 = 1+O(\sqrt\e).$
Moreover, since $\Big(\frac{1+\e}{\sqrt\e}  - \frac{1}{1+\e}
  +4\sqrt\e \Big) =   \frac{1}{\sqrt\e} (1-\sqrt\e + O(\e)),$
we have
\[
\begin{split}
&   \Big(\frac{1+\e}{\sqrt\e}  - \frac{1}{1+\e}
  +4\sqrt\e \Big)
  (1+O(\sqrt\e)) \le
\frac{1}{\sqrt\e}(1+C_1\sqrt\e ),
\end{split}
\]
where the constant $C_1$ is positive and $\e$ is small enough.

Therefore \eqref{pin} becomes
\[
\begin{split}
t & >\frac{2}{\pi} \frac{1}{\sqrt\e}\left[ 1-4\sqrt\e +O(\e) -
\frac{\sqrt\e}{2}(1+O(\sqrt\e) ) \right]
-\frac{2}{\pi}\frac{1}{\sqrt\e} \left[  1+ C_1\sqrt\e +
O(\e)\right]\,\,|z|
\\&
\ge \frac{2}{\pi} \frac{1}{\sqrt\e}\left[ 1-C_2\sqrt\e
  \right]
-\frac{2}{\pi}\frac{1}{\sqrt\e} \left[  1+ C_1\sqrt\e \right]\,|z|.
\end{split}
\]
The latter implies $ |z| > -C_3\sqrt\e   t+ 1-C_4 \sqrt\e, $ which
immediately provides  \eqref{mamma} (note that $(z;t)\in
 B_{1+\e}$ ensures
 $|t|\le C$).   Step  B.2. is finished.

\medskip \noindent \it Proof of Step B.3. \rm
 Now we know that $\ga(1) =(z;t)$ satisfies \eqref{mamma}.
  To better understand the situation, note that,
   as $\e\to 0$, \eqref{mamma} becomes $|z|\ge 1$. Together with the
   shape of the unit ball, which is convex near $|z|=1$, this
    suggests that
  the point $(z;t)$  stays  near the circle $|z|=1$, $t=0$, as $\e$
  approaches $0$. To make this statement quantitative,
   recall also that
our point $\ga(1)$  belongs to the set $ S_{\eta_1}$ by the
biLipschitz property,  $\eta_1\in[(1+\e)^{-1}, 1+\e]$.
 To prove  \it{Step B.3} \rm  write, in a
neighborhood of $|z|=\eta_1$ and $t=0$, the sphere $S_{\eta_1}$ in the form
$
 |z| = \eta_1 \, v(t/\eta_1^2),
$
where $|z|=v(t)$ is the equation of $S_1$ near $|z|=1$, which satisfies
\eqref{toscana}.  Recall
that $\eta_1\in((1+\e)^{-1}, (1+\e))$. Thus we have
the estimate
\[
|z| <(1+\e)\, v\Big( \frac{t}{(1+\e)^2}\Big)\le (1+\e) \left[1 - C
\frac{t^2}{(1+\e)^4}\right]= 1+\e -C\frac{t^2}{(1+\e)^3}
.
\]
This,  together with    \eqref{mamma}, implies
 the estimate in the second line of
(\ref{torino}).
 The proof of Proposition \ref{casimira}
is completed. \endp

\noindent\it Proof of Theorem \ref{bibibi}. \rm
  It suffices to prove it for $I=[0,1]$.
  Write  again $(z;t)=\ga(1)$, where $\g$ satisfy the ODE $\dot\g
= a X(\g)+ b Y(\g)$, $\g(0) = (0,0,0)$ with   $\sqrt{a(s)^2 + b(s)^2 }\in
[(1+\e)^{-1},(1+\e)]$ a.e.  Then
$
z=z(1) = \int_0^1 (a(s), b(s)) ds.
$
Therefore the first estimate of \eqref{torino} (with $s=1$) implies
\[
1-C\sqrt\e \le  \Big| \int_0^1 (a(s), b(s)) ds\Big| \le 1+ C \e.
\]
The proof is finished.  \endp

\section{BiLipschitz image of a horizontal plane}
\setcounter{equation}{0} In this section we prove Theorem B.
This requires the understanding of how the
different quasigeodesics, $(se^{i\phi};0)$ and $(se^{i\theta};0)$,
$s\in \R$, are transformed by $f$, as $\theta, \phi\in [0,2\pi]$.
The key point is in the   following geometric result.

\begin{proposition}
\label{votantonio} Define $\varrho(\theta)=
d\big((1;0),(e^{i\theta}; 0)\big)$, $\theta\in[-\pi, \pi]$. The
function $\r$ is even, smooth on $\left]0,\pi\right[$ and for any
$\la >0 $ there is $C_\la >0$ such that:
\begin{equation}\label{girotondo}
 \r'(\theta)>C_\la \theta^{-1/2} \quad\text{  for all
$\theta\in\left] 0,\pi-\la \right[$.}
\end{equation}
\end{proposition}
 An  immediate consequence of
\eqref{girotondo} is the estimate
\begin{equation} \label{joao} |\r(\theta)-
\r(\phi)|\ge C'_\la  \big|\theta-|\phi|\big|,\quad \forall
\,\theta\in[0,\pi-\la], \phi\in[-\pi+\la ,\pi-\la] .
\end{equation}

The fact that $\r$ has a maximum at $\theta=\pi$ suggests that
 estimate \eqref{joao} no longer holds  for $\la=0$.

We postpone the proof of the Proposition to the second part of the
section.

\medskip \noindent\it Proof of Theorem B. \rm  By
dilation invariance it suffices to prove the statement for  $R=1$.

\smallskip\noindent \it Step 1. Proof of estimate \eqref{expo} for
$|z|=1$.  \rm After a rotation we may assume by Proposition
\ref{casimira}, that $f(1,0,0) = (1+b\e_1; b\e_2 ).$ Observe that
an information on the position of the point $f(-1;0)$
 can be easily extracted. Indeed, write as usual
 $f(1;0)=(\z(1;0),\t(1;0))$.
Formula
 \eqref{vanga} applied  on the interval $(-1,1)$ gives
 \begin{equation}
 \label{rogo}
|\z(-1;0)-(1+b\e_1)| = |\z(-1;0)- \z(1;0)| =
\abs{\int_{-1}^1 c}
=2+b\e_1.
 \end{equation}
 Here we denoted $c(s)= \frac{d}{ds}(\z(s;0))$.
 Moreover, \eqref{torino} gives $
|\z(-1;0)|= 1+b\e_1$, which means $\z(-1;0) = (1+b\e_1)e^{i\psi}$,
for some  $\psi$. Inserting into \eqref{rogo} we get $\psi =
\pi+b\e_2.$ Therefore
 \begin{equation}
 \label{dfgh}
f(-1;0)=\big( (1+b\e_1 )e^{i(\pi+b\e_2) } ;b \e_2\big).
 \end{equation}

In order to prove the required estimate \eqref{expo},
 we will
prove that, after possibly applying the isometry
$(x,y,t)\mapsto(x,-y,-t)$, we can write
\begin{equation}
\label{posizione} f(e^{i\theta};0) = ((1+ b\e_1 )e^{i(\theta +b
\e_3 )};b \e_2 ),\quad\theta\in[-\pi,\pi].
\end{equation}
Note that, by \eqref{rodomonte},  \eqref{posizione}  implies
$d\big(f(e^{i\theta};0), (e^{i\theta};0) \big)\le C\e_4$,   which
is \eqref{expo} when $|z|=1$ and   $A=I$.

Note first
that we already know that (\ref{posizione}) holds for $\theta=0$ and $\theta=\pi$. This
follows from the assumption $f(1;0) = (1+b\e_1; b\e_2 )$ and from
\eqref{dfgh}, which implies
\begin{equation}
\begin{split}
\label{antonietta} d(f(1;0), (1;0)) & \le C\e_3 \qquad \text{and}
\\
  d(f(-1;0),(-1;0)) & =d\big(((1+b\e_1)e^{i(\pi+b\e_2 )};b\e_2
),(-1;0)\big) \le C \e_3 . \end{split}
\end{equation}

We first prove \eqref{posizione} for $\theta\in [0,\frac34\pi]$.
Estimate \eqref{joao} with $\la=\pi/4$ will be used.
 By the results of the previous section
we may write $  f(e^{i\theta};0)= \big( (1+b\e_1 )
e^{i\phi(\theta)}; b \e_2 \big), $
 where
the function $\theta\mapsto\phi(\theta)$ is defined by the last
equality and
 satisfies $\phi(0)=0$. After possibly applying the isometry
$(x,y,t)\mapsto(x,-y,-t)$ we may assume that $\sin
(\phi(\pi/2))>0$, i.e. the second coordinate of $e^{i\phi(\pi/2)}$
is  positive. The biLipschitz property gives
\[
d\big(  f(  e^{i\theta};  0 ), f(1;0)  \big)= (1+b\e) d\big(
(e^{i\theta}; 0), (1;0) \big) = \r(\theta)(1+b\e),
\]
by the  definition of $\r$. By the triangle inequality and the first
line of \eqref{antonietta}, we also have
\[
\begin{split}
 d\big(
 f(e^{i\theta};0), f(1;0)
 \big) & =
 d\big(
 f(e^{i\theta}; 0), (1;0)
 \big) + b \e_3 =
 d\big(
 ( (1+b\e_1)e^{i\phi(\theta)}; b\e_2 ),(1;0)
 \big) +b\e_3
\\&= d\big((e^{i\phi(\theta)}; b\e_2),(1;0) \big) +b\e_3 =
\r(\phi(\theta))+ b\e_3.
\end{split}
\]
Therefore we have proved that $ \r(\phi(\theta))= \r(\theta)+b
\e_3. $ Thus, estimate \eqref{joao} gives $
\big|\theta-|\phi(\theta)| \big|\le C \e_3
 $.
 Since the function
$\phi$  is continuous   and $\phi(\pi/2)>0$, we can drop  the
absolute value:
\begin{equation} \label{ca1}
\big|\theta-\phi(\theta)\big|\le C \e_3, \quad
\theta\in\Big[0,\frac34\pi\Big].
\end{equation}

The same argument works for $\theta\in[-\frac34\pi,0]$
and  estimate \eqref{ca1} also holds in the latter interval.

In order to prove \eqref{posizione} for the  values of $\theta$
near $\pi$ (say $\pi/2\le|\theta|\le\pi$), an analogous argument
can be used,  changing the ``central" point $(1;0)$ with its
opposite $(-1;0)$, whose image's position is narrowed down by the second
line of \eqref{antonietta}. Step 1 is concluded.

\medskip\noindent \it Step 2. Proof of \eqref{expo}
 for
$|z|\le 1$. \rm  We now assume (\ref{posizione}).
   Also, we may assume that
$(z;0)=(r;0)$, $r\in[0,1]$. We know from (\ref{rodomonte}) that
$d(f(r;0),(r;0))\approx
|\zeta(r;0)-r|+|\tau(r;0)+2r\mbox{Im}\,\zeta(r;0)|^{1/2}$, thus we
can estimate the two summands separately. We begin with
$|\zeta(r;0)-r|$. Let $P$ be the point on the segment between $O$
and $\zeta(1;0)$ such that $|P-O|=r$.
 Since, by (\ref{posizione}), the angle
 with vertex $O$ and rays $O\zeta(1;0)$, $O(1;0)$ has amplitude
$b\epsilon_3$ and we have the relations $|(r;0)-O|=r$, we have
$|P-(r;0)|= b\e_3$.
 Consider now the case when
$r\ge1/2$. First we estimate the angle $\alpha$ having vertex in
$O$ and rays $O\zeta(1;0)$, $O\zeta(r;0)$. We claim that
$|\alpha|=b\epsilon_2$.  Indeed, by the Generalized Pythagorean (GPT) Theorem,
\[
|\z(1;0)-\z(r;0)|^2=|\z(1;0)|^2 + |\z(r;0)|^2 -
2 |\z(1;0)|\,\,|\z(r;0)| \cos(\a).
\]
But now, by Theorem \ref{bibibi}, we have
$|\z(1;0)-\z(r;0)|=(1-r)(1+b\e_1)$, $|\z(1;0)|= 1+b\e_1$ and
$|\z(r;0)| = r(1+b\e_1)$. Inserting these estimates into the
previous equation and taking $r\ge 1/2$ into account, we get
$|1-\cos\a|\le C\e_1$, which ensures  $\a=b\e_2.$ Again the GPT
applied to the triangle $O\zeta(r;0)P$ gives
$|P-\zeta(r;0)|=b\epsilon_2$. Therefore,
 the triangle
inequality in the plane gives $|(r;0)-\zeta(r;0)|=b\epsilon_3$.

In the case $r\le1/2$ we proceed much the same way, considering
the triangle $P\zeta(1;0)\zeta(r;0)$ and its angle $\beta$ having
vertex in $\zeta(1;0)$ instead, in order to have the estimate for
$|P-\zeta(r;0)|$.

Finally,  to estimate the second term,
$|\tau(r;0)+2r\mbox{Im}\,\zeta(r;0)|^{1/2} $, observe first that
$|\t(r;0)|\le C\e_2$, if  $r\le 1$. We also know now that $\z(r;
0) =r(1+b\e_1)e^{ib\e_3}$, so that $|\Im \z(r,0)|\le C \e_3$.
Hence $|\tau(r;0)+2r\, \Im\zeta(r;0)|^{1/2}\le C \epsilon_4$. This
ends the proof of Theorem B.
 \endp

\medskip
\noindent\it Proof of Proposition  \ref{votantonio}. \rm Recall
first that   the geodesic balls with center at the origin are
radial in both $|z|$ and $|t|$, i.e. $d_0 (z;t) = d_0(|z|;|t|) $.
The group law gives
\begin{equation}
\label{uuu} \r(\theta) = d\big((1;0),(e^{i\theta};0)\big) =
d_0\big( 2\sin (\theta/2); 4\sin  (\theta/2) \cos(\theta/2) \big).
\end{equation}
 The equation of the upper
half of the sphere $S_r$ is given by  \eqref{eqball} and has the
explicit form $ |z|^2  = \frac{2}{\a^2}(1-\cos\a)r^2,$
 $        t= \frac{2}{\a^2}(\a-\sin\a)r^2$,
 where $          \quad 0<\a<2\pi$. Here we use the coordinate  $|z|^2$  instead of $|z|$,
 in order to make computations easier.
 It is convenient to introduce the functions $G$ and $g$ by the equations:
\begin{equation}
\label{pp1}
 G(\a, r)   =\Big( \frac{2}{\a^2}(1-\cos\a)r^2,
\frac{2}{\a^2}(\a-\sin\a)r^2 \Big) := r^2 g(\a), \quad r>0,\quad
0<\a<2\pi.
\end{equation}
Moreover, the point $(z;t)$ appearing in the right hand side of
\eqref{uuu} satifies  $  |z|^2 =(2\sin (\theta/2))^2 $, $
          t= 4\sin  (\theta/2) \cos(\theta/2) ,$
where $      0<\theta<\pi$.
 Define  the path
\begin{equation}
\label{pp2}   H(\theta) = \Big( 4\sin^2\Big(
\frac{\theta}{2}\Big), 4\sin\Big(  \frac{\theta}{2}\Big) \cos\Big(
\frac{\theta}{2}\Big) \Big) = 2\big(1-\cos\theta,
\sin\theta\big),\quad 0\le \theta\le \pi.
\end{equation}
Observe that $H(\theta)$ describes
the upper half of the circle of radius $2$ centered at $(2,0)$.

By definition of $H$ and $G$, $\rho(\theta)$ is the unique number
with the property that
\begin{equation}
\label{gianfilippo} G(\a, \r(\theta)) = H(\theta)
\end{equation}
for
some $\alpha\in[0,2\pi]$. In fact, (\ref{gianfilippo})  uniquely determines  $(\alpha,\rho)$ as a function of $\theta$. To see  this,
observe that $G(\alpha_1,\rho_1)=G(\alpha_2,\rho_2)$ only when $(\alpha_1,\rho_1)=(\alpha_2,\rho_2)$, otherwise we would have either
two intersecting metric spheres with the same center and
different radii, or a point on a metric sphere whose distance from the center is realized by geodesics
with different values of the parameter $\phi$.

 The proof of Proposition \ref{votantonio}   is articulated as follows.

\medskip
\noindent
 \it{Step 1.} \rm There exist  $C_0, C_1 >0$ such that
$C_1\theta^{1/2}\le \r(\theta)\le C_0\theta^{1/2} $ for any
$\theta$ in $[0,\pi]$.

\noindent
 \it{Step 2.} \rm
 $\r$ is smooth on $\left]0,\pi\right[$ and
$\r'(\theta)$ is strictly positive for any $\theta\in
\left]0,\pi\right[$.

\noindent
 \it{Step 3.} \rm   There exist $\s_0 >0$, $C_0 >0$ such that
$\r'(\theta)\ge C_0\theta^{-1/2} $ for any $\theta\le\s_0 $.

\medskip

\noindent\it Proof of Step 1. \rm  $H$ parametrizes a circle with
speed $2$ and it is easy to verify that $\frac{4}{\pi}\theta\le
|H(\theta)|\le 2 \theta$, for any $\theta\in[0,\pi]$.  On the
other side we have the estimate
\[
  \rho(\theta)^2    \inf_{[0,2\pi]}|g| \le  |
G(\a,\r(\theta))| \le \r(\theta)^2 \sup_{[0, 2\pi]}|g|,
\quad\forall \a\in[0,2\pi],
\]
where $g$ is defined in \eqref{pp1}. The
 required inequalities
follow from the fact that $\displaystyle{0< \inf_{[0,2\pi]}|g| }$
$\displaystyle{< \sup_{[0,2\pi]}|g| <+\infty}$.

\medskip
\noindent\it Proof of Step 2. \rm
 Let $\theta_0\in\left] 0, \pi\right[$. Write
$\r_0=\r(\theta_0)$. Then we have for a suitable $\a_0\in \left]0,
2\pi\right[$  the equation $ G(\a_0, \r_0) = H(\theta_0).
$
The idea is to study the equation
(\ref{gianfilippo})
for
 $\theta$ near a value  $\theta=\theta_0$. We already  know  that there is a unique solution   $(\a(\theta), \r(\theta))$ for  any
 $\theta$ in $[0,\pi]$.   Moreover
  we will show that the function $\r$
 satisfies $\r'(\theta_0)>0$.

In order to apply the inverse function theorem to the  function
$G$, which is smooth near $(\a_0, \r_0)$ we compute
\begin{equation}
\label{irio}
\begin{split}
\p_r G(\a_0 , \r_0 ) & = \frac{4 \r_0}{\a_0^2}(1-\cos\a_0,
\a_0-\sin\a_0)\quad\text {and}
\\
  \p_\a G(\a_0, \r_0) & = -\frac{4 \r_0^2}{\a_0^3} (1-\cos\a_0, \a_0-\sin\a_0)
  + \frac{2 \r_0^2}{\a_0^2}(\sin\a_0, 1-\cos\a_0).
\end{split}
\end{equation}
Then
\[
 \det[\p_\a G(\a_0, \r_0)\quad \p_r G(\a_0, \r_0)]
  = -  \frac{32
 \r_0^3}{\a_0^4}\Big[\sin\Big(\frac{\a_0}{2}\Big) -
 \Big(\frac{\a_0}{2}\Big)\cos \Big(\frac{\a_0}{2}\Big)\Big]\sin\Big( \frac{\a_0}{2}\Big).
\] It is easy to see that the function in the square bracket is
 strictly positive for any $\a_0\in \left]0, 2\pi\right[$.
 Thus, by the inverse function theorem,
  equation \eqref{gianfilippo} can be solved for any $\theta$ near
  $\theta_0$. Denote by $\a(\theta)$, $\r(\theta)$ the solutions.
  The functions $\theta\mapsto \a(\theta)$ and $\r(\theta)$
   are smooth near $\theta_0$.

  In order to get the estimate $\r'(\theta_0)\neq 0$ differentiate
    equation \eqref{gianfilippo}. This gives
$ \p_\a G(\a_0, \r_0) $ $\a'(\theta_0) + \p_r G(\a_0,
\theta_0)\r'(\theta_0) = H'(\theta_0). $
 By  Cramer's rule
\begin{equation}
\label{cramer} \r'(\theta_0) = \frac{\det[\p_\a G(\a_0, \r_0)\quad
H'(\theta_0)] }{\det[\p_\a G(\a_0, \r_0)\quad  \p_r G(\a_0,
\r_0)]}.
\end{equation}
Observe that the second line of \eqref{irio} can be simplified as
follows:
\[
\p_\a G(\a_0, \r_0) = \frac{8 \r_0^2}{\a_0^3} \Big[
\Big(\frac{\a_0}{2}\Big)\cos \Big(\frac{\a_0}{2}\Big)
-\sin\Big(\frac{\a_0}{2}\Big)
 \Big]\Big(\sin\Big(\frac{\a_0}{2}\Big),-\cos \Big(\frac{\a_0}{2}\Big)\Big
 ).
\]
Therefore
\[ \det[\p_\a G(\a_0, r_0) \quad H'(\theta_0)] =\frac{16
\r_0^2 }{\a_0^3}\Big[
 \Big(\frac{\a_0}{2}\Big)\cos \Big(\frac{\a_0}{2}\Big)
 -\sin\Big(\frac{\a_0}{2}\Big)\Big]\sin\Big(\frac{\a_0}{2}+\theta_0\Big),
\]
so that
\begin{equation}
\label{trippa} \r'(\theta_0) =\frac{1}{\r(\theta_0)}\,\, \frac {
\a_0/2 }{\sin(\a_0/ 2)} \sin\Big(\frac{\a_0}{2}+\theta_0\Big).
\end{equation}

Now we are in a position to prove that $\r'\neq 0$. Assume by
contradiction that $\r'(\theta_0)=0$ for some
$\theta_0\in\left]0,\pi \right[$. Then it must be $\frac{\a_0}{2}
+ \theta_0 = \pi $. Equation $G(\a_0, \r_0)=H(\theta_0)$ and the explicit form
of $G$ and $H$ immediately furnish
 \[
  \frac{1-\cos\a_0}{\a_0-\sin\a_0}  =\tan\Big(\frac{\theta_0}{2}\Big)
  = \tan \Big(\frac\pi 2-\frac{ \a_0}{ 4}\Big)
   =  \frac{\cos(\a_0/4)}{\sin(\a_0/4)}
 =\frac{\sin(\a_0/2)}{1-\cos(\a_0 /2)}.
 \]
Observe that   $\sin(\a_0/4)\neq 0$ as $\a_0\in\left
]0,2\pi\right[ $.
 Now
let $\a_0/2= s\in\left]0,\pi\right[$. Then
\[
\frac{2\sin^2 s}{2s -2\sin s\cos s}= \frac{\sin
s}{1-\cos s} .
\]
 But it is easy to see
  that the latter fails
for any $s\in \left] 0, \pi\right[. $  Therefore $\r'\neq 0$ and
Step 2 is accomplished.

\medskip\noindent\it Proof of Step 3. \rm
We are interested in studying equation \eqref{gianfilippo} near
$\theta =0$. If (\ref{gianfilippo}) holds, then
\[
L(\theta):=\frac{1-\cos\theta}{\sin\theta}=
\frac{1-\cos\a}{\a-\sin\a}:=R(\a).
\]
Note that $L(\theta)= \tan\left(\theta/2\right)$ increases on
$[0,\pi]$, from $L(0)=0$ to $L(\pi-)=+\infty$.  One readily
verifies that $R$ is a strictly monotone function decreasing from
$R(0^+)=+\infty$ to $R(2\pi)=0$. Hence $\alpha$ is a monotonically
decreasing function of $\theta$, $\alpha(0)=2\pi$, $\alpha(\pi)=0$.
Keeping into account that $\theta = 2\arctan R(\a)$,
 a calculation shows that (i)
$\frac{d\alpha}{d\theta}<0$ for all $\theta\in \left]0,\pi\right]$; (ii)
$\frac{d\alpha}{d\theta}= -\frac{\pi}{2\pi-\alpha}(1+o(1))$, as
$\theta\to0$ $(\Leftrightarrow\,\a\to 2\pi$); (iii)
$\frac{d\alpha}{d\theta}=-\frac{3}{2}(1+o(1))$, as $\theta\to\pi $
 ($\Leftrightarrow \a\to 0$).
As a consequence,
there are $C_1$ and $C_2>0$ so that,
 \begin{equation}
 \label{pentola}
 C_1\sqrt\theta\le\frac{2\pi-\alpha}{2}\le  C_2\sqrt\theta.
 \end{equation}
 Finally we go back to  \eqref{trippa}. For
$\theta$ close to $0$ (which means $\frac{\a}{2}$ close to $\pi$),
we have the estimates $\sin\frac\a 2\le \Big( \pi -\frac \a
2\Big)$   and $\sin \Big(\frac\a 2+\theta\Big)\ge
C_3\Big(\pi-\frac\a 2-\theta\Big)$.   
Therefore,
\[
\r'(\theta)\ge  \frac{1}{\r(\theta)}
 \frac{C}{\pi-\frac\a 2}\left(\pi-\frac{\alpha}{2}-\theta\right).
\]
The latter, together with (\ref{pentola}) and the estimate
$\r(\theta)\simeq\theta^{1/2}$, as $\theta\to 0$, concludes the
proof of Proposition \ref{votantonio}.
\endp

\section{Image of points outside the  plane $t=0$}
\setcounter{equation}{0}

Let $f$ be a biLipschitz map  as in Theorem B. By Theorem B we
know how  the plane $t=0$ transforms: for any $R>0$ there is a
suitable $A\in O(2)$ such that,
\begin{equation}
\label{pisello} d((f(z;0), (A z ; 0)) \le   C\e_4 R , \quad |z|\le
R.
\end{equation}

In the current section, in order to study where points outside the
plane $\{t=0\}$ are mapped, we will make a systematic use of the
following family of geodesics $s\mapsto \g(s)= (x(s), y(s),
t(s))$, where
\begin{equation}\begin{split}
\label{modi}
  & x(s) =
q  \cos\a \Big(1-\cos\Big(\frac sq \Big) \Big) -q\sin\a
\sin\Big(\frac sq\Big),
\\
 &y(s) =  q\sin\a  \Big(1-\cos\Big(\frac sq \Big) \Big)
 +q\cos\a  \sin\Big(\frac sq\Big),
\\&
 t(s) =  2 q^2 \Big(\frac sq - \pi -\sin\Big(\frac sq\Big)\Big).
\end{split}\end{equation}
The parameter $q$ is positive, while $\a\in [0, 2\pi]$. Note that
for $0\le s/q\le 2\pi$, $d(\gamma(s),O)\approx q$.
 The path $\g$ is
a unit speed geodesics with lifetime $2\pi q$. It can be obtained
from \eqref{pietro} with a translation and by changing $\phi$ with
$1/q$. Moreover
\begin{equation}
\label{cerri} \g(0) = (0,0, -2q^2\pi),\quad \g(\pi q)= (2q\cos \a,
2q\sin\a, 0), \text{ and }  \g(2\pi q) = (0,0, 2q^2\pi).
\end{equation}
The distance $\dist\big((0;2\pi q^2),\{ t=0 \}\big)$ is realized
by any point of the form $(2qe^{i\theta };0)$
 and its value is
\[
\dist\big((0;2\pi q^2)),\{ t=0 \}\big)= d\big( (0;2\pi q^2),
(2qe^{i\theta};0) \big) =\pi q, \quad\forall\theta\in[0,2\pi].
\]
All the points of  the circle $(2qe^{i\theta};0)$,
$\theta\in[0,2\pi]$ are  ``projections" of  $(0; 2\pi q^2)$ on the
plane $\{ t=0\}$. This is the reason why statement 1. in
Proposition \ref{sesto} below is false for $s=0$.

 By means of an
accurate analysis of the mentioned geodesics, we will obtain the
following quantitative result, whose technical proof will be given in Subsection \ref{tecc}.

\begin{proposition}\label{sesto}
There exist  universal constants $\s_0>0$ and $C_0>0$
 such that, for any $\s<\s_0$
 the following statement holds.
 For any  $q\in\left ]0,\infty\right[$  consider the unit speed geodesic $\gamma$
 of total lifetime $2\pi q$ such that $\gamma(0)= (0,0,-2q^2\pi)$ and
 $\gamma(\pi q)= (2q, 0,0):=Q $. Take any number $s$ with
 \begin{equation}
\label{gianbattista} \s^{1/8}\le s/ q \le \pi - \s^{1/16}
 \end{equation}
 and denote $P=\gamma(s) $. Then:

1. The closure of the  ball $B(P, d(P,Q) )$  touches the plane
$t=0$ only in $Q$.

2. Let $\Pi_1(z;t)=z$   be the vertical projection on $\{t=0\}$.
The enlarged ball $ B\big(P, (1+\s)d(P,Q)\big)$ satisfies the
following property.
\begin{equation}
\label{peste}
\Omega_1= \Pi_1 \Big( B\big(P, (1+\s)d(P,Q)\big) \cap
\{t>0\}\Big)
 \subset \{ z\in\C :
|z-2q|\le C_0  q \s^{1/4} \}.
\end{equation}
\end{proposition}

\begin{figure} [ht]
 \centerline{\includegraphics{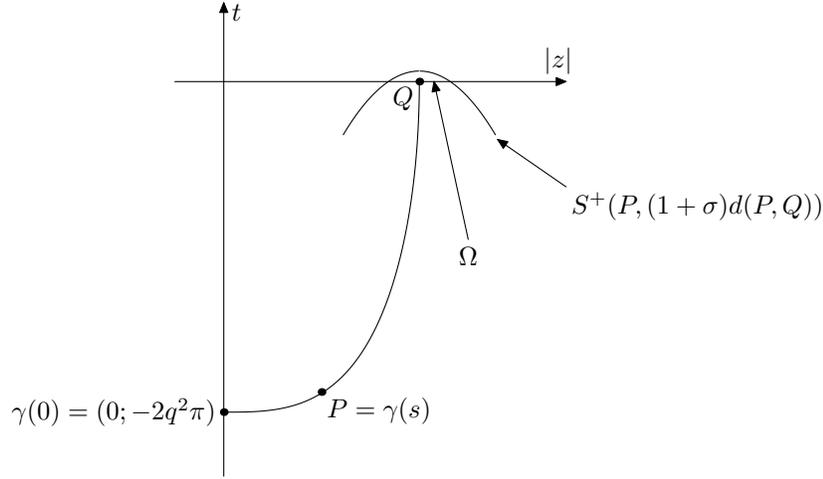} } \caption{ A bidimensional qualitative 
representation
 of inclusion \eqref{peste}.}
\end{figure}

\subsection{Points on the $t-$axis}
Next  we analyze the position of points of the form $f(0;t)$. Our
result is the following
\begin{theorem}\label{toto}
There are $\e_0>0$ and $C>0$ such that, if $\e<\e_0$, $(z;t)\mapsto f(z;t) =
(\z(z;t); \t(z;t))$ is $(1+\e)-$BiLipschitz on $\H$ and satisfies
$f(0)=0$, then
\begin{equation}
\label{out} |\z (0,0,t)|\le C\e_4 |t|^{1/2}
 \quad\text{and}\quad
\big|\left|t\right| - |\t(0,0,t)|\big|\le C\e_4 |t|,\quad \forall \quad t\in\R.
\end{equation}
\end{theorem}
The proof of Theorem \ref{toto} involves only values of $f$ on the
set $ K:=\{t=0\}\bigcup \{$ $t-$axis$\}$. This does not contain
information enough to determine the sign of $\t(0,0,t)$ (the map
  $f(x,y,t)=(x,y,-t)$, which is far from being an isometry in $\H$, is an 
isometry while
  restricted to $K$).
 In order to determine the sign of $\t(0,0,t)$, we need to take into account  
values of
$f$ outside the set  $K$.  This is done in Proposition \eqref{roux} below.

\begin{proposition}\label{roux}
If $f$ is $(1+\e)-$biLipschitz, $\e\le\e_0$, $f(0)=0$ and $\det(A)
= +1$ in \eqref{pisello} for some $R>0$,  then    
\[
\frac{\t(0,0,t)}{t}>0 \quad \text{for all $t\in\R$, $t\neq 0$}.                 \]
\end{proposition}
  Before giving the proof of   Theorem \ref{toto}   and Proposition
 \ref{roux}, observe that, putting  the mentioned statements together, we  
 immediately get  the proof of Theorem C.

\bigskip \noindent \it Proof of Theorem C.  
 \rm Let $f$ be a 
$(1+\e)-$biLipschitz self map of $\H$ such that $f(0)=0$.   By Proposition \ref{roux}, after
possibly applying a rigid motion $(x,y,t)\mapsto (x,-y,-t)$,
we may delete the absolute value in the second inequality in \eqref{out}  Thus,
\[ 
              |\z (0,0,t)|\le C\e_4
|t|^{1/2}
 \quad\text{and}\quad
\big|t - \t(0,0,t)\big|\le C\e_4 |t|,\quad \forall \quad t\in\R,
 \] 
which implies 
\begin{equation}
\label{pisellonio2}
d(f(0;t), (0;t)) \le C \e_5 |t|^{1/2},\quad \forall t\in\R,
\end{equation}
as desired. \endp

\medskip

\noindent \it Proof of Theorem \ref{toto}. \rm
  Since the statement
is dilation invariant, we prove it for the point $(0;2\pi)$. By
Theorem B with $R=2$, we may assume that
\eqref{pisello} holds with $A=I$ for $R= 2$. Write $f(0;2\pi) =
(\xi,0,\t)$, $\xi>0$ (the proof which follows can be easily
modified to cover the general case $f(0,0,2\pi)=(\xi,\eta,\t)$).
 Note that $(0,0,2\pi) = \g(2\pi )$, where $\g$ is
one among the geodesics  in \eqref{modi}, with $q=1$, see
\eqref{cerri}. Moreover, we know that
 the distance of the point
$(0;2\pi  )$  from the plane $t=0$ is realized by all the points
of the form $(2  e^{i\theta};0) $ and its value is $\pi $.

  The
idea   of the proof is the following. We will choose
$\theta=\pi/2$ and  $\theta=\frac 32\pi$. We will show that, by the biLipschitz
property, the point $ (\xi, 0,\t) $  has
 distance $\pi  +o(\text{power of }\e )$ from  both the points $(0, 2, 0)$
 and $(0, -2, 0)$. These information, together with the the one about the distance
  of  $(\xi, 0, \t)$ from the origin,
  $d(f(0;2\pi), (0;0))=
(1+b\e)d((0;2\pi), (0;0))$, will give a rigid estimate of the
position of the point $(\xi, 0,\t).$

Take $\theta=\frac{\pi}{2}$. By the triangle inequality and the
biLipschitz property we have
\begin{equation} \label{abaco}
\begin{split}
d\big( (\xi,0,\t), (0,2, 0) \big)& =d\big(f(0;2\pi), (0,2,0)\big)
\\&=
d\big( f(0;2\pi), f(0,2, 0) \big) + b d\big(f(0,2, 0), (0,2, 0)
\big)
\\&
=\pi (1+ b\e) +b \e_4 =\pi+b\e_4,
\end{split}
\end{equation}
where we used \eqref{pisello}, which holds for $R=2$.
 The same computation for the
opposite point $(0,-2 , 0)$ shows that
\begin{equation}
\label{abaco2} d\big( (\xi,0,\t), (0,- 2  , 0) \big)= \pi
(1+b\e_4).
\end{equation}
Write again \eqref{abaco} and \eqref{abaco2} using the group law
and recalling that the distance from the origin satisfies
$d((0;0), (z;t)):=d_0(z;t) = d_0(|z|; |t|)$ for any $(z;t)\in\H$.
This gives
\begin{equation}
\begin{split}
\label{dinto}
 d_0\left(\sqrt{\xi^2
+4 },\t-4\xi  \right) =\pi  (1+b\e_4)=
 d_0\left(\sqrt{\xi^2
+4 },\t+4\xi   \right).
\end{split}
\end{equation}
Denote now  
$
\r=\sqrt{\xi^2 +4 } ,$ $ \t_-=\t-4\xi$ and $\t_+=\t+ 4\xi  .$ The
equivalences in \eqref{dinto} can be written as
\begin{equation}
\label{dintoski} d_0(\r;\t_+)= \pi (1+b \e_4)\quad\text{and}\quad
d_0(\r;\t_-)= \pi (1+b\e_4 ).
\end{equation}

Next we prove that $|\t|\ge \t_0 $ for some small but absolute
constant $\t_0>0$, uniformly for small $ \e$. In order to get this
property   we add to \eqref{dinto} (or the equivalent
\eqref{dintoski}) the third information given by the biLipschitz
property
\begin{equation}
\label{dintello} d_0(\xi,0,\t)=d_0\big(f(0;2\pi) \big) =(1+ b \e)
d_0(0;2\pi ) = (1+ b \e) \, \pi \sqrt 2
\end{equation}
(see  \eqref{assicuro},   for the last equality). By
\eqref{dinto}, since the ball $B(0,r)$ is contained in the
cylinder $ \{|z|< r \}$ for all $r>0$, we have
$
\xi^2 + 4  \le \pi^2 (1+C\e_4).$ Since  $\pi^2- 4 <2\pi$, this
gives  for small $\e$ the estimate
\begin{equation}
\label{emloski} \xi^2\le 2\pi , \quad \text{that is}\quad \xi\le
\sqrt{2\pi} .
\end{equation}
It is immediate to see  (note that $\pi \sqrt 2 >\sqrt{2\pi}$)
that \eqref{emloski} and \eqref{dintello} together imply
\begin{equation}
\label{doc3} |\t|\ge \t_0 ,
\end{equation}
for some absolute constant $\t_0$.

Now, given $C_0>0$ introduce the  ring domain 
\[
A_\e:=B\left(0,\pi  (1+ C_0 \e_4 )\right)\setminus B\left(0,\pi
(1-C_0\e_4 )\right) , 
\] and let $ A_\e^+= A_\e\cap\{ t>0\} $,  $
A_\e^-= A_\e\cap\{ t<0\} $. By \eqref{dintoski}, we may choose an
absolute constant $C_0 > 0$ such that $(\r; \t_+), (\r; \t_-)\in
A_\e$. Note that \eqref{doc3} says that  if $\e$ is small enough,
it can not happen
 that $(\r;\t_+)\in A_\e^+$ and $(\r;\t_-)\in
A_\e^-$. Thus it should be $(\r;\t_+), (\r;\t_-)\in A_\e^+$ or
$(\r;\t_+), (\r;\t_-)\in A_\e^-$.

It is not difficult   to check the following claim by means of the properties of the unit ball described in Section \ref{iniziale}.

\medskip\noindent\it Claim. \rm If $\t_0>0$ is given, then there exists $\s_0$ and $C_0>0$
 such that, for any
$\s<\s_0$, given any pair of point $(\r;\t_-), (\r;\t_+)\in
B(0,1+\s)\setminus B(0,1-\s) $,  $\t_-\le\t_+$ with $\t_+\ge\t_0$
(or $\t_-\le -\t_0$), then $\t_+-\t_-\le C_0 \s$.
 \rm

\medskip

Rescaling  the claim  (with $\s=C_0\e_4  $) from the
unit radius to the radius $\pi  $,  we get the estimate
$
\t_+ - \t_-\le C \e_4, $ which by the definition of $\t_+, \t_-$
gives $4\xi \le C\e_4$, that is $\xi\le C  \e_4$.
  This ends the proof of the first
inequality in \eqref{out}.

In order to prove the second inequality in \eqref{out}, recall
that, by \eqref{dintoski},
\[
\left( \sqrt{\xi^2 +4} ,\t-4\xi  \right) \in S \left(0,\pi  (1
+b\e_4)\right) \quad\text{and }\xi\le C\e_4.
\]
 Inserting these information into  equation \eqref{ornette}
of the sphere   of radius $\pi(1+b\e_4)$,
 we get
\begin{equation}\label{erfro}
\left|\frac{\t +b\e_4}{\pi^2 (1+b \e_4)^2}\right| = u \left(
\frac{\sqrt{4+b\e_3}}
 { \pi (1+b \e_4)}
\right).
\end{equation}
Assume first  that the quantity inside the absolute value is
positive. Recall that,  by \eqref{assono}, $u (w) = \frac 2\pi +
O(w- \frac {2}{\pi} )^2$ as $w\rightarrow {\frac 2\pi}$. After a
short manipulation \eqref{erfro} becomes
$
\t=2\pi +b\e_4  $, which is the required estimate.
 If instead the number in the absolute value in the left hand side of
\eqref{erfro} is negative, then we  get $\t=-2\pi+b\e_4$.
Ultimately,  the second estimate of \eqref{out} holds and the
proof of the theorem is finished.
\endp

\bigskip
\noindent\it Proof of Proposition \ref{roux}. \rm   By dilation invariance it suffices to prove the proposition for $R=2$, i.e. 
assume  that \eqref{pisello} holds
with $A = I_2$, the $2\times 2$ identity matrix, and $R=2$.
 We prove that
the point $P= (0,0,- 2\pi)$ goes into a point whose $t-$coordinate $\t(0;-2\pi)$  satisfies  $\t(0;-2\pi)<0$.  By continuity, since  the second inequality of \eqref{out} ensures 
 that $\t(0;t)\neq 0$, for all $t\neq 0$, this will be enough to prove the proposition.

In the proof of this proposition, which is qualitative, $o(1)$
always denotes (scalar or vector) functions such that $|o(1)|\le
C\e^k$ for some absolute but unimportant positive constants $C$
and $ k$, which may change at each occurrence.

By  \eqref{out} we know that it should be
$
\z(0,0,- 2\pi) =   o(1) $ and $|\t(0,0,-2\pi)|= 2\pi +o(1).$
Assume by contradiction that  $\t(0,0, -2\pi) = +2\pi +o(1)$.
Consider the geodesic \eqref{modi} with $\a=0$ and $q=1$, which
has the form
\begin{equation}
\label{giugio} \g(s)=(1-\cos s, \sin s, 2(s-\pi-\sin s)).
\end{equation}
Note that $P=\g(0).$ Write $Q=\g(\pi) = (2,0,0)$ and note also
that $\g(2\pi) = -P = (0,0,2\pi)$. Take now the intermediate point
$M =(\g(\pi/2))=(1+i;-(\pi +2))$. Our assumption $\t(0;-2\pi) =
+2\pi +o(1)$ implies also $\t(0;-(\pi+2))= +\pi+2 +o(1)$.
Moreover, \eqref{out} gives also  $\z(0;-(\pi+2))=o(1)$. Our
knowledge on global quasigeodesics, applied to the quasigeodesic
$\la(s):= f\big( \frac{1+i}{\sqrt 2}s; -(\pi+2)\big)$, $s\in\R$
(note that $\la(0)=f(0;-(\pi+2))$ and $\la(\sqrt 2)= f(M)$), tells
us that it should be
\begin{equation}
\label{bbv} f(M)=((1+i)e^{i\theta};\pi +2)+o(1),
\end{equation} for
some $\theta\in[0,2\pi]$. Furthermore,  we may also assert that,
by \eqref{pisello}, with  the matrix $A=I_2$, $f(Q)=Q+o(1)$.
  Then, the triangle inequality and the biLipschitz
property give
\begin{equation}
\label{bbc} d(f(M), Q)= d(f(M), f(Q))+o(1)
=d(M,Q)+o(1)=\frac{\pi}{2}+o(1).
\end{equation}
To get some information on $\theta$, we use   Proposition
\ref{sesto}.   Indeed,
since both \eqref{bbv} and \eqref{bbc} hold,
it must be
\begin{equation}\label{mirabilia}
f(M) = \g\left(\frac 32\pi\right)+o(1)=(1,-1,\pi+2)+o(1),
\end{equation}
where $\gamma$ is defined in \eqref{giugio}.  To check
\eqref{mirabilia}, consider the geodesic $\gamma$ restricted to
$[\pi/2,2\pi+\pi/2]$. Since $\gamma(\pi)=Q$,
$f(M)+o(1)\in\{t=\pi+2\}$,
 $\gamma(3\pi/2)\in\{t=\pi+2\}$, 
$f(M) \in B(Q,\pi/2+o(1))$, Proposition \ref{sesto} says that
$d(f(M),\gamma(3/2\pi))=o(1)$, hence (\ref{mirabilia}) holds.

 Finally use the biLipschitz property $ d\big(f(M),
f(1,1,0)\big)= d\big(M, (1,1,0)\big)+o(1)  $, that is
$d\big((1,-1, \pi+2 ), (1,1,0) \big)= d\big((1,1,-(\pi+2)),
(1,1,0)\big) +o(1) $. Translating in term of the distance from the
origin $d_0$,
\begin{equation}
\label{erge} d_0\big(0,-2,\pi-2  \big)= d_0\big( 0,0,-(\pi+2)
\big)+o(1)\quad \Rightarrow d_0\big( 2;\pi-2
\big)=\sqrt{\pi(\pi+2)}+o(1).
\end{equation}
We have concluded that the point $(2;\pi-2)$ has distance from the
origin $\sqrt{\pi(\pi+2)}+o(1)$.
 The latter number  is greater than $4$ if $o(1)$ is small enough.
 But the ball of radius $4$ contains the rectangle
 $\big[0, \frac {8}{\pi}\big]\times \big[0, \frac{16}{\pi}\big]$ (see Section
 \ref{iniziale}). The point  $(2;\pi-2)$ is strictly inside the mentioned
 rectangle. This is in  contradiction with \eqref{erge}. \endp

\subsection{Image of points outside the $t-$axis.}
\begin{theorem}[Theorem D]  
There are $\e_0$ and $C_0$ absolute constant such that, if $f$ is
$(1+\e)-$biLipschitz, $\e\le\e_0$,  $f(0)=0$ and   $ A=I_2$ in
Theorem B  at the scale  $R>0$, then
\[
d\big(f(z;t),(z;t)\big)\le C\e_{11} R,
\]
for any $(z;t)$ s.t. $[|z|^2+|t|]^{1/2}\le C_0 R$.
\end{theorem}

\noindent\it Proof. \rm We prove the statement for $R=1$. Consider
a point $P=(z;t)$, outside the set $\{ t=0\}
\bigcup\{t-\text{axis}\}$. To locate quantitatively the position
of $f(P)$, we will use Proposition \ref{sesto}. Therefore it is
convenient to think the point $P$ in the form $P=\ga(s)$, where
$\ga$ is the geodesic in \eqref{modi}, for some $q>0$ and $0<s<\pi
q$. We may also choose $\a=0$. The choice  $R=1$ means $q\approx
d_O(z;t)\le C_0$, $C_0$ absolute.

 Roughly
speaking, if the point is near $z=0$ or near the $t-$axis, we will
get the required estimate by means of the previous results and the triangle inequality. If instead the point is
``far" from  $\{ t=0\} \bigcup\{t-\text{axis}\}$, then we will
invoke Proposition \ref{sesto}.

To be more precise, we will distinguish the following cases:

\medskip
\noindent  \it Case A.1: \rm $0\le q\le\e_{11} $  ($(z;t)$  close to  the origin).

\noindent\it Case A.2: \rm $  \e_{11}\le q\le C_0 $ and $
0\le\frac{s}{q}\le \e_9
 $ ($(z;t)$    close to the $t-$axis).

\noindent \it Case A.3: \rm $\e_{11}\le q\le C_0 $ and $  \pi -\e_{10}\le
\frac{s}{q}\le \pi $ ($(z;t)$ close to the plane $\{ t=0\}$).

\noindent \it Case B: \rm  $\e_{11}\le q\le C_0 $ and $\e_9\le \frac
sq\le\pi-\e_{10}$ ($(z;t)$   far from $\{t=0\}\bigcup
t-\text{axis}$).

 We discuss
first the cases $A$, which all will be treated by the triangle
inequality.

\smallskip\noindent
\it Case A.1. \rm Recall that $f(0)=0$ and $\g(s)=(z;t)$. The
triangle inequality gives
\[
\begin{split}
d\big(f(z;t), (z;t)\big) & \le d\big(f(z;t), (0;0)\big)+
d\big((0;0), (z;t)\big)
\\&
  \le (2+\e) d\big((z;t),(0,0)\big)\le C\e_{11} .
\end{split}
\]

\smallskip
\noindent\it Case A.2. \rm We have $0\le\frac  sq\le\e_9$. Then
\[
\begin{split}
d\big(f(z;t),(z;t) \big) & \le d\big(f(z;t),f(0;t)\big)+
d\big(f(0;t), (0;t) \big)+d\big((0;t), (z;t) \big)
\\&
 \le (2+\e) d\big((z,t), (0;t))\big)+C\e_5\le C|z|+C\e_5,
\end{split}
\]
where we used biLipschitz property, triangle inequality and
\eqref{pisellonio2}.    Moreover, since $\frac s q\le \e_9$,
\eqref{modi} gives, for small $\e$, $|z|\le
C\e_9 $. Therefore the right-hand side  can be estimated by
$C\e_9$ which is clearly smaller than $C\e_{11}$.

\smallskip\noindent
\it Case A.3. \rm Use the triangle inequality and
\eqref{stimantonio}.
\[
\begin{split}
d\big(f(z;t), (z;t)\big)&\le d\big(f(z;t), f(z;0)\big)+
d\big(f(z;0), (z;0)\big) + d\big((z;0), (z;t)\big)
\\&
  \le (2+\e) d\big((z;0), (z;t)\big)+C\e_3 \le(2+\e)\sqrt{\pi
  t}+C\e_3.
\end{split}
\]
Since $\pi\ge \frac sq \ge\pi -\e_{10}$,
 we have, by \eqref{modi}, if
$\e$ is small enough, $|t|\le \e_{10}$. Therefore the last line
can be estimated by $C \e_{11}$.

\smallskip\noindent
\it Case B. \rm Write again $(z;t)=\g(s)$ and, as usual
$f=(\z;\t)$. The key point is to show that, since, by hypothesis,
\eqref{pisello} holds with $A=I$ and $R=1$, then
\begin{equation}
\z(z;t)=z+b\e_8 \quad\text{and}\quad\t(z;t)=t+b\e_4 ,\label{ans}
\end{equation}
for any $(z;t)$ such that $|z|^2+|t|\le C_0$ and Case B holds.

 To prove \eqref{ans}, recall that we know that $f(0;t)=
(b\e_4; t+b\e_4) $, by Theorem \ref{toto}. Then, by our result on
the image of a horizontal plane, Theorem B,
\[
f(0;t)^{-1}\cdot f(z;t)= \big( z(1+b\e_1)e^{i\b}; b\e_2\big),
\]
for some $\b\in [0,2\pi]$. Therefore, writing $f(z;t)= f(0;t)\cdot
 \big(  z(1+b\e_1)e^{i\b}; b\e_2\big) $, it turns out that
$\t(z;t)= t+b\e_4$. This is the second equality in \eqref{ans}.

In order to get the first one, we need to locate $\z(z;t)$ with
the help of Proposition \ref{sesto}. This will provide information
on the angle $\beta$. Recall first that, if $Q=\gamma(\pi q)$,
then  $d(f(Q), Q)\le C\e_4$.
 The triangle
inequality,  the
   biLipschitz assumption and the (already proved) second equation of \eqref{ans} give
\begin{equation}
\label{opto}
\begin{split}
 d\big( (\z(z;t);t) ,  Q\big)
 &
 \le d\big( (\z(z;t);t) ,f(z;t) \big) +
  d\big(f(z;t), f(Q) \big) +d\big(f(Q),Q \big)
\\&
 \le C \e_5 +(1+\e)d((z;t),Q)+C\e_4 \le  d((z;t), Q)+ C\e_5.
 \end{split}
\end{equation}
To write \eqref{opto} in a  form which is more suitable for the application of
Proposition \ref{sesto}, recall that (Case B) we are assuming
$q\ge\e_{11}$ and $\pi-\frac{s}{q}\ge\e_{10}$. Then, since
$(z;t)=\ga(s)$ and $Q=(\ga (\pi q)$,  we have
\[
d((z;t), Q)=\pi q-s =q\Big(\pi-\frac sq\Big)\ge\e_{11}\e_{10}\ge
\e_6.
\]
 Then $\e_5=\e_6^2\le \e_6 d((z;t), Q)$.
Thus \eqref{opto} takes the more dilation invariant form
\begin{equation}
\label{retia} d\big( (\z(z,t);t) ,  Q\big)\le
d((z;t),Q)\{1+C\e_6\}.
\end{equation}
Looking at \eqref{retia} and recalling $\z(z;t) =
ze^{i\beta}+b\e_1$, we get by triangle inequality that
\[
\begin{split}
&Q\in B((ze^{i\beta};t), d((z;t), Q) (1+C\e_6)),\ \mbox{i.e.}
\\&
\quad \mathcal{R}_{-\b}Q\in B((z;t) , d((z;t), Q) (1+C\e_6)),
\end{split}
\]
where $\mathcal{R}_{-\b}(w;s)=(e^{-i\b}w;s)$.
This, together with the assumption  of Case B, which provides
\eqref{gianbattista}   with $\s= C\e_6$,
enables us to apply
  apply \eqref{peste} with $\s= C \e_6$, which reads
$
  |R_{-\b}Q- Q|\le C q \e_6^{1/4} = C q\e_8.
$
Dividing both members by $q$ gives $|e^{i\beta}-1|\le
C\epsilon_8$, hence $ \big|\z(z;t)-z
\big|\le   C \e_8.
$
Thus    \eqref{ans} is proved.

Finally, write $f(z;t)$ in the form given by \eqref{ans}. Then
\[
d\big(f(z;t), (z;t)\big)= d_0\left((-z; -t)\cdot \big( z+b\e_8;
t+b\e_4 \big)\right)\le C  \e_9\le C\e_{11},
\]
as desired. The proof of the theorem is concluded. \endp

\subsection{Proof of Proposition \ref{sesto}} \label{tecc}
Observe first that  statement 1. follows from the proof of Lemma
\ref{convesso}.  In order to prove
statement 2., note that, by dilation invariance, we may  choose
$q=1$. Letting $\a=0$ and $q=1$ in the geodesic \eqref{modi} gives
 \[ P= \g(s) =
\left( \big(1-\cos s \big), \sin s
 , 2  (s-  \pi -\sin s) \right).
\]
It is  $d(P,Q)= \pi - s$. Put $R:= (1+\s)d(P,Q) = (1+\s)(\pi -s)$.

Before proving (\ref{peste}), we show that $\Omega_1=\Omega$,
where $\Omega$ is the seemingly smaller set $$ \Omega:=  B\big(P,
(1+\s)d(P,Q)\big)\cap\{ t=0 \}. $$ The equality $\Omega_1=\Omega$
holds if   
the surface $S_{\mathrm{up}}=S(P,R)\cap\{t>0\}$ can be
viewed as the graph of a function, $(z;t)\in S_{\rm{up}}$ iff $t=t(z)$,
and this is true if and only if the ``equator'' $E$ of $S(P,R)$
lies in $\{t\le0\}$. By equator, we mean the set $E=P\cdot
(S(O,R)\cap\{t=0\})$, the set of the points in $S(P,R)$ which have
parameter $\phi=0$ (see Section 2). The equator of $S(O,R)$ has
parametrization $(R\cos\alpha,R\sin\alpha,0)$, $0\le\alpha<2\pi$.
After a left translation by $P$, we see that the $t$-coordinate of
a point in $E$ has equation
\begin{eqnarray*}
t(s,\alpha)&=&2(s-\pi-\sin s)+2R\sin s\cos\alpha-2R\sin\alpha(1-\cos s)\cr
&=&2(s-\pi-\sin s)+4R\sin(s/2)\cos(s/2+\alpha)\cr
&\le&2(s-\pi-\sin s)+4R\sin(s/2)=k(s).
\end{eqnarray*}
We want then to show that $k(s)\le0$ when $\s^{1/8}\le s \le \pi - \s^{1/16}$.
 Passing to the coordinate to $u=\pi-s$ and replacing 
$R$ by its explicit expression,
the inequality holds if
\begin{equation}
\label{barbagatto}
0<u+\sin u-2u(1+\sigma)\cos(u/2)=h(u),\ \mbox{if}\ \sigma^{1/16}\le u\le\pi.
\end{equation}
A Taylor expansion shows that, for $u=\sigma^{1/16}$, (\ref{barbagatto}) becomes
$$
2\sigma^{1/16}-\frac{1}{4}\sigma^{1/16+1/8}(1+o(1))
<2\sigma^{1/16}-\frac{1}{6}\sigma^{3/16}(1+o(1)),
$$
which is true if $\sigma<\sigma_0$ is small enough.
For the other values of $u$, we take a derivative in (\ref{barbagatto}),
$$
h^\prime(u)=2\cos^2(u/2)-2\cos(u/2)+{u}\sin(u/2)-b\sigma=g(u)-b\sigma.
$$
Observe that $g(\sigma^{1/16})=1/4\cdot\sigma^{1/8}(1+b\sigma)$ and that
$$
g^\prime(u)=2\sin(u/2)-2\sin(u/2)\cos(u/2)+\frac{u}{2}\cos(u/2)\ge0\ \mbox{when}\ 0\le u\le\pi.
$$
Hence,
$
h'(u)=g(u)-b\sigma\ge g(\sigma^{1/16})-b\sigma>0,
$
if $\sigma^{1/16}\le u\le\pi$. This shows that $k(s)\le0$, hence that $\Omega=\Omega_1$.

\smallskip

 We now return to the proof of (\ref{peste}) for $\Omega_1=\Omega$.
The generic point  $A$ of the upper half sphere
 $S^+(0,R)$ has coordinates
 \begin{equation}\label{agramante}
A(\alpha,\phi)=\left(2\frac{\sin(\phi R/2)}{\phi}\cos\a,
2\frac{\sin(\phi R/2)}{\phi}\sin\a,
 2\frac{\phi R -\sin(\phi R)}{\phi^2}\right)
, \quad
\end{equation}
where $   |\a| \le 2\pi,$ and $   0\le\phi R\le 2\pi. $
Denote by
$(x,y,t)$ the coordinates of the point $P\cdot A :=(x,y,t)\in
S^+(P, R)$. Then
\begin{equation}
\begin{split}
\label{giornale} x&= (1-\cos s )+ 2\frac{\sin(\phi
R/2)}{\phi}\cos\a, \qquad y
  = \sin s +
2\frac{\sin(\phi R/2)}{\phi}\sin\a
\\
t&= 2 ( s-\pi)  +\frac {2R}{\phi } - 2\Big\{ \sin s
 +\frac{\sin (\phi R)}{ \phi^2} \Big\}
+\frac{8}{\phi}\sin \Big(\frac s2 \Big) \sin\Big(\frac{\phi
R}{2}\Big)
 \cos
 \Big(\a+\frac{s}{2}\Big).
\end{split}
\end{equation}
Note that letting $\a=-\frac{s}{2}$, $\phi=1$ and $\s=0$,
 i.e.  $R=\pi-s$, we have   $t=0$ for all $s\in(0, \pi)$, as expected.

To prove the proposition, take a small $\s$, fix $s$ satisfying
\eqref{gianbattista} and   consider the function $t=t(\a, \phi)$
defined in the last line of \eqref{giornale}. We prove the
following two statements.

\medskip \it\noindent Step 1. \rm  For small $\s$, the
$t-$coordinate corresponding to
  $\a=-\frac{s}{2}$ and  $\phi =1$ is positive.

\smallskip\noindent
\it
Step 2. \rm For any point of the form
\begin{equation}\label{jedi}
\mu(\psi)=(\a,\phi)
=\left(-\frac{s}{2} +\s^{1/4}\cos \psi, 1+ \s^{1/4}\sin \psi
\right),
\end{equation}
 it is  $t<0$  for any   $\psi\in
  [0,2\pi]$.

\smallskip
Once the described steps are proved, we will show that they ensure
the proof of the proposition.

\medskip \noindent \it Proof of Step  1. \rm Put $\a=-\frac s2$, $\phi=1$
 and $R=(1+\s)(\pi-s)$ into the third equation of \eqref{giornale}.
 After some simplifications  and a Taylor
expansion near $\s=0$, we get
\[
\begin{split}
t & =2\s(\pi -s) - 2\sin s -2 \sin (s-\s(\pi -s))+8\sin \Big(\frac
s2\Big)  \cos\Big(\frac s2 -\frac{\s}{2}(\pi - s)\Big)
\\&
  = 2\s(\pi - s)\Big( 1+\cos s + 2\sin^2\Big( \frac  s 2
  \Big)\Big) + o(\s(\pi -s))
 = 4\s(\pi -s)+   (\pi -s)o(\s),
\end{split}
\]
 as $\s\to 0$. Then $t>0$, if $\s$ is smaller than an absolute constant
 $\s_0$.

\medskip
\noindent \it Proof of Step 2. \rm  We introduce the more
comfortable variables $x$, $\b$, $\la$,
\[
x=\pi- s,\qquad \phi=   1  + \b,\qquad \a=  - s/2 +\la .
\]
Then $\phi R =(1+\b)(1+\s)x$. Put also
$
\phi R-x:=\d = (\b+\s+\s\b)x. $ Then, starting again from the last
line of \eqref{giornale}, we get
\[
\begin{split}
\frac{t(1+\b)^2}{2}
 &=
 (1+\b)(\s-\b) x -(1+\b)^2 \sin x  -\sin (x+\d)
 \\&\qquad+
 4(1+\b)\cos\Big(\frac x 2\Big)\sin\Big(\frac x 2+\frac \d
 2\Big)\cos \la.
\end{split}
\]
Expanding the right hand side for $\s$ and $\d$ near $0$, we get
\[
\begin{split}
\frac{t(1+\b)^2}{2}
 &=
O(\s) -\b x -\b^2 x -\sin x -2\b \sin x -\b^2\sin x
\\&
 \qquad -\sin (x) -\cos (x) \d +\sin x \frac{\d^2}{2}+ O(\d^3)
 \\&\qquad+
 4\cos\Big(\frac x 2\Big)\Big\{  (1+\b)\cos\la\Big( \sin
 \Big(  \frac x 2\Big)+\cos  \Big(  \frac x 2\Big)\frac \d 2-\sin
  \Big(  \frac x 2\Big)\frac{\d^2}{8} + O(\d^3)\Big)\Big\}_1 .
\end{split}
\]
 Another Taylor expansion at the second order in $\b,\la,\d$, gives
\[
\begin{split}
\{ \dots \}_1 =\sin  \Big(  \frac x 2\Big)+ \b\sin  \Big(  \frac x
2\Big) +\cos  \Big(  \frac x 2\Big)\frac\d 2 +\cos  \Big(  \frac x
2\Big) \b\frac\d 2 -\sin  \Big(  \frac x
2\Big)\frac{\d^2}{8}-\frac{\la^2}{2} \sin  \Big(  \frac x 2\Big)
+R_1,
\end{split}
\]
with the estimate  $R_1\le C_0(|\b|^3+|\la|^3 +|\d|^3)$, as soon
as $\b +\la+\d\le \s_0 $. $C_0$ and $\s_0$ are absolute constants.
Observe in particular  that all these expansions are   uniform in
the variable $x\in [0,\pi]$.
 Now recall that $\d =
(\b+\s+\s\b) x=\b x +O(\s)$. Then we can write $\b x$ instead of
$\d$ and $\b^2 x^2$ instead of $\d^2$, making an error of $O(\s)$.
Then
\[
\begin{split}
\frac{t(1+\b)^2}{2}
 &= -\b^2\Big\{  \sin x-x\cos x-\frac{x^2}{4}\sin
 x\Big\}-\la^2\sin x + R_2,
\end{split}
\]
with  $R_2\le C_0(\s+|\b|^3 +|\la|^3)  $, as before, if $\b
+\la+\d\le \s_0 $.

To conclude the argument, note that the function  $ x\mapsto \sin
x-x\cos x-\frac{x^2}{4}\sin
 x$ is increasing on   $(0,\pi)$ (it has positive derivative).
 Therefore, since it behaves as $C_1x^3$, $C_1>0$, near $0$, it turns out that
 $ \sin x-x\cos x-\frac{x^2}{4}\sin
 x\ge C_1 x^3$,  for any $0\le x\le \pi$.   Then
\[
\frac{t(1+\b)^2}{2}\le -\b^2 C_1 x^3 -\la^2 \sin x + C_2\{
\s+|\b|^3 +|\la|^3\} ,
\]
provided $\s,\b,\la$ are small enough.  Now hypothesis
\eqref{gianbattista}, in term of our variable $x$, becomes $x\ge
\s^{1/16} $ and  $x\le \pi-\s^{1/8}$, so that it is also $x^3\ge
\s^{3/16} $ and  $\sin x\ge C_3\s^{1/8}$. Write $\b = \s^{1/4}\cos
\psi$, $\la= \s^{1/4} \sin\psi$. Then
\[
\frac{t(1+\b)^2}{2}\le- C_1\s^{11/16 }\cos^2 \psi -
C_3\s^{5/8}\sin^2 \psi + C_2 \s^{3/4}.
\]
Now  if  $\s$ is small enough with respect to the absolute
constants $C_j$'s, we have proved that $t<0$.  This ends the proof
of \it{Step 2}. \rm

\smallskip

Now, if the Carnot sphere were convex,
then Steps 1 and
 2 would give almost immediately inclusion \eqref{peste}. This is not the case, but we will show that all
of the interesting action takes place in the convex part of the
ball's boundary.

We claim that $A(\mu(\psi))\subset S^+(O,R)\cap\partial
B_{co}(O,R)$, where $B_{co}$ denotes the convex envelope, in the
Euclidean sense, of $B(O,R)$. By Lemma \ref{convesso}, this
amounts to showing that $|A_1(\mu(\psi))|\ge\frac{2R}{\pi}$,
where  we decomposed $A=(A_1;A_2)$. Clearly $
A_1(\alpha,\phi)=2\frac{\sin(\phi R/2)}{\phi}e^{i\alpha} $, where
 $(\phi, \a)$ are of the form  \eqref{jedi}. This implies
   $\phi\in[1-\s^{1/4}, 1+\s^{1/4}]$. Moreover,  $R= (\pi-s)(1+\s)$
   and
  $s$ satisfies \eqref{gianbattista}. Then we have $\phi R/2\in
\left]0,\pi/2\right[$. Thus, the elementary  inequality
$\sin(x)\ge\frac{2}{\pi}x$,   $x\in[0,\pi/2]$, provides the required lower  estimate on $|A_1|$.

Let $\mu^*=(\mu^*_1; \mu^*_2 ): =P\cdot A(\mu)$. After
a translation, the fact that the curve $A(\mu )$ lies in the
convex part of $S^+(O,R)$ implies that
$\Omega^\prime=B_{co}(P,R)\cap\{t=0\}$ is a convex set contained
inside the curve $\mu^*_1$. {\it A fortiori,}
$\Omega\subset\Omega^\prime$ is contained inside $\mu^*_1$.

To finish the proof, we have to show the inclusion in
(\ref{peste}). Since $\Omega$ is contained inside $\mu^*_1$,
it suffices to prove that
\begin{equation}
\label{zerlina} |\mu^*_1(\psi)-2|\le C_0\sigma^{1/4},\
\psi\in[0,2\pi],
\end{equation}
for some absolute $C_0>0$. After a translation, this amounts to
showing that $|A_1(\mu(\psi))-A_1(-\frac{s}{2},1)|\le
C_0\sigma^{1/4}$. The latter follows from the definition
\eqref{jedi} of $\mu$ and the elementary estimate
$|D_{Euc}A_1(\alpha,\phi)|\le C R\le C$ for the Euclidean
derivative's norm of $A_1$.
\endp

\section{Approximation of derivatives}\setcounter{equation}{0}
\label{mirabilandia}
We prove here the following theorem:
\begin{theorem}\label{bmo}
There are universal constants $C>0$ and $\e_0>0$ such that the following
statement holds. Let $f$ be $(1+\e)-$biLipschitz with $\e\le \e_0$,
$f(0)=0$.  Assume that for some $r>0, k\in\N,$
\[
d\big( f(z;t), (z;t)\big)\le C\e_{k} r,\quad (z;t)\in B(0,r).
\]
 Denote by $J f$
the Jacobian matrix of $f$ in the sense of Pansu. Then,
\begin{equation}
\label{casello} \int_{B(0, C^{-1}r)}\|Jf(x,y,t)-I \|dxdydt\le C
{\e_{k+1}} r^4.
\end{equation}
\label{ivana}
\end{theorem}
In the above  estimate, $I=\begin{pmatrix}
 1 & 0 \\
  0 & 1
\end{pmatrix}$ and $\|\cdot\|$ denotes the norm of a  $2\times 2$
matrix.

By the elementary properties
 of Pansu derivative, see Section
\ref{iniziale}, if a given map $f = (\xi,\eta,\t)$
 is Pansu-differentiable at a point   $P $  and $Jf =\begin{pmatrix}
   \a & \b\\
   \g & \d \
 \end{pmatrix}$, then   we have
\begin{equation}
\label{glu}
 \frac{d}{ds}\xi\big(e^{sX}(P)\big)\Big|_{s=0}=X\xi(P)=\Big(
Df \big(P\big)(1,0,0) \big) \Big)_1= \a(P).
\end{equation}
Here we used the notation  $(x,y,t)_1 = x$  for the first
component and $s\mapsto e^{sX}(P)$ denotes the integral curve of
$X$ emanating from $P$ at $s=0$. An analogous formula holds for
$Y$.

\bigskip
\noindent \it Proof of Theorem \ref{ivana}. \rm
 It is not restrictive to choose   $r=1$.
First we show  that if $f$ is $(1+\e)-$biLipschitz, then the
Jacobian matrix $Jf $ satisfies
\begin{equation}
\label{glo} Jf(P)^T Jf(P) = I_2+b\e, \qquad \text{for a.e.
$P\in\H$,}
\end{equation}
where the entries of  $b$ are bounded by an absolute constant.
In particular its diagonal  components, whose sum plays the role
of the divergence that appears in John's proof, satisfy the estimate
\begin{equation}
\label{glin} |\a(P)| \le 1+  C_1\e,\qquad |\d(P)| \le 1+  C_1\e, \qquad
\text{for a.e. $P$.}
\end{equation}
To prove \eqref{glo}, recall that for a.e. $P$ there are $\a,
\b,\g,\d$
s.t.
\[
\big(\a u+ \b v, \g u+\d v, (\a\d-\b\g)w\big)=Df(P)(u,v,w)  =
 \lim_{\s\to 0}\delta_{1/\s}\big\{f(P)^{-1}\cdot
f\big( P\cdot \d_\s(u,v,w) \big)\big\}.
\]
Then, letting as usual $d_0$ to indicate the distance from the
origin
\[
\begin{split}
d_0 \big( \a u+ \b v, \g u+\d v, (\a\d-\b\g)w \big)
&= \lim_{\s\to 0} d_0\left(
\delta_{1/\s}\big\{f(P)^{-1}\cdot f\big( P\cdot \d_\s(u,v,w)
\big)\big\} \right)
\\&
= (1+ b \e)d_0(u,v,w),
\end{split}
\]
by the biLipschitz property. Taking $w=0$ we get, at any
differentiability point $P$,
\begin{equation}
\label{commendatore}
\left| \binom{u}{v}\right|(1+b\e) = d_0(u,v,0)(1+b\e)= d_0 \big(
\a u+ \b v, \g u+\d v, 0 \big) = \left| Jf(P)\binom{u}{v}\right|,
\end{equation}
for all $(u,v)\in \R^2$. Equality \eqref{glo} then follows by
simple considerations of linear algebra in the Euclidean plane.
 Finally,
 \eqref{glin} follows immediately from
\eqref{glo}.

Our next  task   is to  follow John's argument, starting from a
one dimensional estimate and integrating it by Fubini's Theorem.

Consider the set $\Omega_1=\{e^{sX}(0,y,t): |y|,|t|,|s|\le 1 \}$.
Here and in the following we denote by $e^{sX}(0,y,t)$ the
integral curve of $X$ starting at $s=0$ from the point $(0,y,t)$.
The map $(s,y,t)\mapsto   e^{sX}(0,y,t)= (s, y, t+2ys)$ is volume
preserving. By  the change of variable formula and Fubini Theorem,
we get the formula
\begin{equation}
\label{figghiu} \int_{\Omega_1} g = \int_{|y|, |t| \le 1} dy
dt\int_{-1}^1 ds \,g(e^{sX}(0,y,t)),
\end{equation}
for any function $g$.  By Pansu's Theorem, there is $\Sigma\subset
[-1,1]\times[-1,1]$ of full 2-dimensional measure such that, given
any $(y,t)\in \Sigma$, the map  $f$ is Pansu differentiable at the
point $e^{sX}(0,y,t)$ for a.e. $s\in[-1,1]$.

Introduce  the function  
\[
\Phi(x,y,t)=f(x,y,t)^{-1}\cdot (x,y,t):= \big( u(x,y,t), v(x,y,t),
w(x,y,t)\big).
\]
Note immediately that $u(x,y,t)=x-\xi(x,y,t)$,
$v(x,y,t)=y-\eta(x,y,t)$. 
In spite of the fact that $\Phi$  may be neither Lipschitz, nor Pansu differentiable,
  see Remark \ref{remo}, we can  
define at any $P$ where $f$ is Pansu differentiable, the $2\times 2$ matrix
\begin{equation*}
 \wt J\Phi(P) :=  I_2 -J f (P),
\end{equation*} or, letting $\wt J\Phi 
=  \begin{pmatrix}
\wt\a    &\wt\b \\
   \wt\gamma & \wt\d
 \end{pmatrix}$,
\[
 \begin{pmatrix}
\wt\a(P)    &\wt\b(P) \\
   \wt\gamma(P) & \wt\d(P)
 \end{pmatrix}= \begin{pmatrix}
 1-\a(P)    & -\b(P) \\
   -\gamma(P) & 1 -\d(P)
 \end{pmatrix}.
\]
By the formula for $Jf$ in \eqref{posticino}, given $(y,t)\in\Sigma$, we have for almost any $s\in[-1,1]$,
\[
\begin{split}
\frac{d}{ds}u(e^{sX}(0,y,t) ) & = \frac{d}{ds} \Big( s - \xi
(e^{sX}(0,y,t)) \Big) = 
\text{by \eqref{glu}}
\\&=1-\a(e^{sX}(0,y,t)) \ge 1-|\a(e^{sX}(0,y,t))| \ge 
-C_1\e,
\end{split}
\]
by estimate \eqref{glin}.
  Therefore $
\frac{d}{ds}  u(e^{sX}(0,y,t))+C_1\e\ge 0$.

Now we are ready to integrate: take $(y,t)\in\Sigma$. Then
\[
\begin{split}
\int_{-1}^1| \wt\a(e^{sX}(0,y,t))|ds &=
\int_{-1}^1   \abs{\frac{d}{ds}  u(e^{sX}(0,y,t))} ds
\\
& \le
\int_{-1}^1\Big\{ \Big|  \frac{d}{ds}  u(e^{sX}(0,y,t))+C_1\e
\Big|+C_1\e\Big\}ds
\\&
 = \int_{-1}^1\Big\{   \frac{d}{ds}  u(e^{sX}(0,y,t))+2 C_1\e\Big\}ds
\\& =u(e^X(0,y,t))- u (e^{-X} (0,y,t)) +4C_1\e.
\end{split}
\]
By hypothesis  $d(f(P), P)\le C\e_{k}$. Since
$(f(P)^{-1}\cdot P)_1= u(P)$,
 the first two terms
can be estimated by $C \e_{k}$. Thus
\[
\int_{-1}^1 |\wt\a (e^{sX}(0,y,t))| ds \le C\e_{k},\quad\text{for a.e.
$ (y,t)\in[-1,1]^2$}
\]
Integrating over $[-1,1]^2$, see \eqref{figghiu},
$
\int_{\Omega_1}|\wt\a|\le C\e_{k}.
$
The same
argument can be used with the field $Y$ instead of $X$. Then we
get, for a suitable $\Omega\subset\Omega_1$,
\[
\int_{\Omega}|\wt\a |+|\wt\d |\le C\e_k.
\]
 The estimate of the trace of   $\wt J\Phi$
is accomplished.

The remaining part of the proof can be concluded exactly as in
John's paper (see \cite[p. 407]{J}). We just sketch it.  Recall
that $\wt J\Phi = I-Jf$. Put
$
U = \wt J\Phi+(\wt J\Phi)^T $ and $  V=(\wt J\Phi)^T \wt J\Phi.
$ Thus
$
\|V-U \|=\| (\wt Jf)^T Jf -I  \|\le C\e,
$
by \eqref{glo}. Ultimately
\[
\int_{\Omega}\|\wt J\Phi\|^2=\int_\Omega \|V\|\le
\int_\Omega|\mathrm{tr}(V)| \le C\int_\Omega|\mathrm{tr}(U)|\le
C\e_k.
\]
It now  suffices to apply H\"older's inequality $\int_\Omega\|\wt J\Phi\| \le
C\big(\int_\Omega\|\wt J\Phi\|^2\big)^{1/2} $.
The proof is concluded. \endp

\begin{remark}\label{remo}
 Consider $f(x,y,t)=(2x,y,2t)$, $(x,y,t)\in\H$. The function $f$ has Lipschitz 
constant  $2$. The corresponding  
map $\Phi(P)=f(P)^{-1}\cdot P$ has the form
\[
 \Phi(x,y,t)= f(x,y,t)^{-1}\cdot(x,y,t)=(-x,0,-t+2xy).
\]
Let $x\ne0$.
Testing the Lipschitz condition for $\Phi$ for the points $P=(x,y,t),$ $Q=(x,y+\delta,t-2x\delta)$ as $\delta\to0$, we see that the Lipschitz constant of $\Phi$ at $(x,y,t)$ can not be finite.
Moreover, $\Phi$   is not Pansu 
differentiable at $P.$
\end{remark}

Theorem \ref{bmo} says that $Jf$ belongs to $BMO(\H)$. In the Euclidean case, the John-Nirenberg inequality \cite{JN} allowed John to deduce
a local exponential integrability result for the Jacobian of a biLipschitz map.
In the context of the Heisenberg group the same conclusion holds, due to the far-reaching
generalization of the John-Nirenberg inequality due to Buckley \cite{Bu}.
\begin{corollary}\label{ircocervo}
There exist
 constants $\epsilon_0,C>0$ such that, if $\epsilon\le\epsilon_0$, $f$ is $(1+\epsilon)$-biLipschitz on $\H$ and $B$ is a ball
in $\H$, then
\begin{equation}
\label{mangiafuoco}
\frac{1}{\mathcal{L}(B)}\int_{B}\exp\left(\frac{\|Jf(Q)-(Jf)_B\|}{C\epsilon_{12}}\right)dQ\le2.
\end{equation}
\end{corollary}
  In (\ref{mangiafuoco}), $(Jf)_B$ is the average of $Jf$ on $B$ and the constant $2$ on the right hand side could be replaced by any
constant $\lambda>1$, changing the value of $C$.

\section{Examples}\label{olimpia}
In this section we  discuss some examples.

\begin{example} \rm
We show here   that, in Theorems C and D, the powers of $\e$ on
the right hand side of the inequalities can not be improved to be
$\epsilon^1$, but have to be at least $\e^{1/2}$. Consider the
dilation 
\begin{equation}
 \label{dilao}
f(z,t)=\delta_{1+\e}(z;t)=((1+\e) z;(1+\e)^2 t )
\end{equation} which is $(1+\e)$-biLipschitz in
reason of the homogeneity of the distance function. It is obvious
that $d(f(0;1),(0;1))=c\sqrt{\e}d((0;0),(0;1))$ for some explicit
absolute $c>0$. This shows that the estimate of Theorem C can not
hold with $\e^1$ in the right-hand side.

Next,  we consider again  the function $f$ in \eqref{dilao} and  we show that for any isometry $\Gamma$ of $\H$ there is a
point $P$ such that $d(P,O)\le\sqrt{\pi}$ and
$d(f(P),\Gamma(P))\ge c\sqrt{\epsilon}$ for some absolute $c>0$.
By the proof of the isometries' classification in the Appendix,
any isometry $\Gamma$ of $\H$ can be written as
$\Gamma=L_{(w,s)}\circ {\mathcal R}_\theta \circ J^m$,   for some
$\theta\in{\mathbb R}$, $(w,s)\in\H$ and $m\in\{0,1\}$. We assume
$m=0$, the other case  being
similar. We have $\Gamma(0;0)=(w;s)$ and $\Gamma(0;1)=(w;s+1)$,
hence $A:=d(f(0;0),\Gamma(0;0))=d_0(w;s)$ and
$B:=d(f(0;1),(\Gamma(0;1)))=d_0(w,s-(2\e+\e^2))$.   From the
geodesic equations, (\ref{pietro}) or (\ref{eqball}), we deduce
that, for fixed $w$,
 $\max\{A,B\}$ is minimized when $s=\frac{1}{2}(2\e+\e^2)$, hence $A=B$.
Keeping $s=\frac{1}{2}(2\e+\e^2)$ fixed, it is easy to see that
$A$ is minimized when $|w|=\sqrt{2s/\pi}$ and $A=\sqrt{\pi
s/2}\sim\sqrt{\pi\e/2}$, for small $\e$.
\end{example}

Next we briefly describes two procedures for producing contact
maps, devised respectively
 by  Kor\'anyi and  Reimann \cite{KR1} and by Capogna and Tang \cite{CT}.
Recall   that a $C^1$ diffeomorphism $f:\H\to\H$ is contact if 
its differential sends  horizontal vectors to horizontal vectors. 
 Before introducing the procedures we observe the following
 standard fact. Let $f:\H\to\H$
   be a $C^1$
 map and assume that $f$ is contact.
    Denote by   $Jf$
    the invariant components of its Jacobian, see formula
    \eqref{posticino}. Then, letting
 \begin{equation}
 \label{pepe}
L =\sup_{p\in\H}\|Jf(P)\|,
 \end{equation}
 the map $f$ is $L-$Lipschitz.
Here $\|\cdot\|$ is the operator norm of the matrix, acting on
Euclidean ${\mathbb R}^2$.
We omit the  standard proof (but see \cite[Theorem 3.2]{AM}, where a more general statement is proved).

\begin{example}
[Kor\'anyi and  Reimann type maps.] \rm In \cite{KR1}, Kor\'anyi
and  Reimann show how to produce quasiconformal maps as flows of a
suitable vector field.

Consider a function $p:\H\to \R$, say   $C^2-$smooth.
Define the vector field
\begin{equation}
\label{campo} v=-\frac 1 4 (Yp)X +\frac 1 4 (Xp)Y +pT,
\end{equation}
Denote by $f_s(P)$, $P\in \H$, the solution of the Cauchy problem
$
\frac{d}{ds}f_s(P)=v(f_s(P)),$ $ f_0(P)=P.$ It is known that such
a vector field generates a contact flow. The differential of  the map
$f_s$ at $P\in\H$ sends $\HH_P$ into $\HH_{f(P)}$. See
\cite[Theorem 5, p.331]{KR1}.

It is not difficult to check, by a slight modification of the argument in \cite{KR1},
 that a condition on $p$ which ensures
the biLipschitz property is an estimate of the form
\begin{equation}
\label{179}
\sup_{\R^3}\{|X^2p|+|Y^2p|+|XYp|+|YXp|\}= C_0<\infty.
\end{equation}
In that hypothesis,   for all $s\in\R$ the map $f_s$ is
biLipschitz and the biLipschitz constant is controlled by
$L=e^{C_1|s|}$.

Note that, in order to obtain an estimate on the Lipschitz
constant, we assume  \eqref{179}, which is slightly stronger than
the one in \cite{KR1}, which involves only a bound on
$\sup|Z^2p|$, $Z=X-iY$.
\end{example}

\begin{example}[Maps which preserve vertical lines.] \rm
We follow \cite{CT}, \cite{BHT}. Consider a nonsingular contact
map $f:\H\to\H$ of the form $f(z;t)=(\z;\t)=(\z(z);\t(z;t))$.
 We say
that $f$ is {\it vlp, vertical lines preserving}. If $f$ is vlp,
then $Jf(P)$ coincides with the Euclidean Jacobian of $\z$.  If $f$
as above is $C^2$, a standard calculation shows that $f$ is vlp if
and only if $\mbox{det}J\z (P)=\gamma$ is constant, $\gamma\ne0$,
and $\t=B(z)+\gamma t$, where $B$ can be recovered from its
Euclidean differential, $1/2dB=(xdy-ydx)-\gamma(udv-vdu)$, where
$x+iy=z$ and $u+iv=\z$.

By (\ref{pepe}),
 we have that the biLipschitz norm of the vlp map
$f=(\z;\t)$ equals the {\it Euclidean} biLipschitz norm of the map
$\z:{\mathbb R}^2\to{\mathbb R}^2$.  Observe that the dilation
$\delta_{1+\epsilon}$ considered above is vlp. Other interesting
examples of vlp maps arise when we consider $\z$ to be one of the
spiral-like plane maps in \cite{GM}. By lifting their maps to $\H$
we obtain, for $k\le0$, the vlp maps: $$ S_k(z,t)=(ze^{ik\log|z|},
t-k|z|^2). $$ By the results in \cite{GM}, $S_k$ is
$\alpha$-biLipschits, with
$\alpha=\frac{|k|+\sqrt{|k|+4}}{2}=1+\frac{5}{8}|k|+o(k)$. The
image of the plane $\{t=0\}$ under $S_k$ is the cone $\{(w,s):\
s=|k||w|^2\}$, hence this class of examples dos not say anything
new on the power of $\e$ in Corollary \ref{ammira}.
\end{example}

\section*{Appendix: the case $\e=0$}
\renewcommand{\thetheorem}{A.\arabic{theorem}}
\renewcommand{\theequation}{A.\arabic{equation}}
\setcounter{equation}{0}  \setcounter{theorem}{0}

We show here that any isometry of $\H$ which fixes the origin has
the form $f(z;t)=J^m{\mathcal R}_\theta$ for some
$\theta\in{\mathbb R}$ and $m\in\{0,1\}$, see Section
\ref{iniziale}. There are at least two proofs in the literature.
The first is by noting that isometries are $1$-Quasiconformal maps
and that the latter are described in \cite{KR1} and \cite{C1}.
 The second consists in analyzing the geometry of the group
 $\H$ at the level of its Lie algebra  \cite{Ki}.

 The simple proof we provide below relies on properties of the distance $d$ alone. We were
 interested in finding such a proof to have a clue
at how to investigate biLipschitz mappings of $\H$ from a purely metric point of view.

Let $f:\H\rightarrow\H$ an isometry such that  $f(0)=0$.
Consider the geodesic
$\gamma(s)=(s,0,0)$, $s\in\R$. This is a globally minimizing
 geodesic and it is sent by the
isometry $f$ to another globally minimizing geodesic, which has
the form  $f(s,0,0) = (sv; 0 )$,  $v_1^2+v_2^2=1.$ Up to a
rotation we may choose  $v=(1,0)$. In other words we may assume
that   $f$ is the identity on the line $y=t=0$. The same argument
shows that the image of the plane $\{t=0\}$ is the plane itself.

Next, look at the set  $S = \{  (x,y,0): x^2 +y^2 = 1 \}$.
Since rotations are isometries and $(-1,0;0),\ O, \ (1,0;0)$ are collinear,  $f(S)= S$, $f(1,0,0)= (1,0,0)$
  and   $f(-1,0,0) = (-1, 0,0)$.
Up to a transformation of the form
$(x,y,t)\mapsto(x,-y, -t)$ we can assume that
 $S^+ := S\cap\{y>0\}$ is mapped onto itself.
 We claim that  $f$ is the identity on $S$. Suppose there is a
 point
$(e^{i\theta_0};0)$,   $\theta_0\in\left ]0,\pi\right[$ which is not
sent onto itself by $f$, say $f((e^{i\theta_0};0)) =(e^{i\theta_1};0)$
with $\theta_1>\theta_0$ (the opposite case can be treated in the
same way). Inductively $(e^{i\theta_{n+1}};0) = f((e^{i\theta_n};0))$.
Since the map is one-to-one, the sequence of angles is strictly
increasing:
 $\theta_0<\theta_1<\cdots <\theta_n <\cdots$. Either there is $\theta_k$ such that $\theta_k<\pi$ and
 $\theta_{k+1}>\pi$, but this contradicts the fact that
 $f(S^+)\subset S^+$,  or the sequence $(\theta_n)_{n\ge 0}$
 is infinite. Since $f$ is an isometry,
\[
d((e^{i\theta_1};0),( e^{i\theta_0};0)) = d((e^{i\theta_2};0),
(e^{i\theta_1};0))=\cdots = d((e^{i\theta_n};0), (e^{i\theta_{n-1}};0))=\cdots
\]
 Now $\theta_n\rightarrow
\bar\theta\le\pi$, which implies $d((e^{i\theta_n};0),( e^{i\theta_{n-1}};0))\to 0$, as $n\to\infty$.
But this contradicts the fact that
 $ d((e^{i\theta_n};0),
(e^{i\theta_{n-1}};0))$ is the same for all $n$.

From the above it follows that, up to a composition with a map 
$J$,  
 the map
  $f$, when restricted to the plane $t=0$, is the identity.

Next consider the plane  $t=\bar t$, where $\bar t> 0$. This
plane is sent in a left translate of the plane $t=0$.
 But the only left translates of $t=0$ which do not intersect  the plane
 $t=0$ itself (this would violate the injectivity of $f$)
  have the form $t=$constant.
 Therefore  $f(\{t=\bar t\}) = \{ t=\wt t \}$, for a suitable  $\wt
t\neq 0$.

Now we claim that  $f(0,0,\bar t) = (0,0,\tilde
 t)$. This follows as  $(0,0,\bar
 t)$ is the unique point of the plane  $t=\bar t$ which can be connected
 through geodesics
 lying in the plane $t=\bar t$ to any other
 point $(z;\bar t)$, $z\in\C$, the same happens to its image. Thus,
  $f(0,0,\bar t) = (0,0,\tilde t)$.
Formula \eqref{assicuro} tells also that
it must be $\tilde t =\pm \bar t$.

Assume first that
  $\wt t = \bar t$ (the opposite case will be discussed later).
  The image of the global geodesic $(s,0,\bar t)$, $s\in\R$,
is of the form  $f(s,0,\bar t) = (sv; \bar t)$,
 where $|v|=1$. To recognize that  $v=(1,0)$, observe that
\[
d((s,0,\bar t), (s,0,0)) =  d(f(s,0,\bar t), f(s,0,0) ) = d
((vs;\bar t),(s,0,0))\quad \forall\, s\in\R.
\]
After a left translation (write as usual $d_0$ for the distance from the
origin),
\[ d_0(0,0,\bar t) = d_0 ((-1+v_1)s, v_2 s, \bar t +
2 v_2 s^2) \quad \forall\, s\in\R.
\]
But this can hold only if $v_2=0$ and $v_1=1$  (otherwise the point
 $((-1+v_1)s, v_2 s, \bar t + v_2 s^2)$ would go to infinity, as  $s\rightarrow\infty$).

Then, if  $f(0,0,\bar t) = (0,0, \bar t)$,   $f$ is the
identity on $t=\bar t$. The same argument can be repeated at any
quote $t=t^*$, $t^*\in\R$ and the proof is finished.

Finally,
consider the case  $f(0,0,\bar t) =
(0,0,-\bar t) $. Arguing as before, write $f(s,0,\bar t)= (vs;-\bar
t)$, $s\in\R$. Then
\[
\begin{split}
d((s,0,\bar t), (s,0,0)) & =
d(f(s,0,\bar t), f(s,0,0))=d((vs;-\bar t),(s,0,0) )
\\&
  =d_0\big((v_1-1)s, v_2s, -\bar t + 2 v_2 s^2\big),\quad s\in\R,
\end{split}
\]
which implies $v_1=1$ and $v_2=0$. Now we discover that in this case $f$ cannot be an isometry.
Without loss of generality suppose $\bar t=1$ and choose
  $(z; \bar t)=(1,0,1) $. This gives
\[
d((1,0,1), (1,1,0)) = d(f(1,0,1), f(1,1,0)) =
d((1,0,-1), (1,1,0)),
\]
which implies $d_0(0,1,1)=d_0(0,1,3)$, a false equality.
 \endp

\small

\noindent Nicola Arcozzi: Dipartimento di Matematica,\\
 Universit\`a di Bologna.\\
 Piazza di Porta San Donato, 5,\\
 40127  Bologna, Italy.\\
E mail: \tt arcozzi@dm.unibo.it\rm

\medskip
 \noindent Daniele Morbidelli: Dipartimento di Matematica,\\
 Universit\`a di Bologna.\\
 Piazza di Porta San Donato, 5,\\
 40127  Bologna, Italy.\\
E mail: \tt morbidel@dm.unibo.it
\end{document}